\documentclass[11pt]{article}
\usepackage{amsfonts}
\usepackage{color}
\usepackage{amsmath,amssymb,comment}
\newtheorem{theo}{Theorem}[section]
\newtheorem{lem}[theo]{Lemma}

\newtheorem{defi}[theo]{Definition}

\newcommand{\mysection}[1]{\section{#1} \setcounter{equation}{0}}
\newcommand{\proof}{{\sc Proof.} \quad}
\newcommand{\proofc}{{\sc Proof} \ }
\newcommand{\be}{\begin{equation} \label}
\newcommand{\ee}{\end{equation}}
\newcommand{\bea}{\begin{eqnarray}\label}
\newcommand{\eea}{\end{eqnarray}}
\newcommand{\bas}{\begin{eqnarray*}}
\newcommand{\eas}{\end{eqnarray*}}
\newcommand{\bit}{\begin{itemize}}
\newcommand{\eit}{\end{itemize}}
\newcommand{\qed}{\hfill$\Box$ \vskip.2cm}
\newcommand{\nn}{\nonumber}
\newcommand{\R}{\mathbb{R}}
\newcommand{\N}{\mathbb{N}}
\newcommand{\pO}{\partial\Omega}

\newcommand{\eps}{\varepsilon}

\newcommand{\supp}{{\rm supp} \, }

\newcommand{\wto}{\rightharpoonup}
\newcommand{\wsto}{\stackrel{\star}{\rightharpoonup}}
\newcommand{\hra}{\hookrightarrow}
\newcommand{\io}{\int_\Omega}
\newcommand{\na}{\nabla}

\newcommand{\al}{\alpha}
\newcommand{\lam}{\lambda}

\newcommand{\sig}{\sigma}
\newcommand{\pa}{\partial}
\newcommand{\bom}{\overline{\Omega}}
\newcommand{\Om}{\Omega}

\newcommand{\wh}{\widehat}

\newcommand{\hs}{\hspace*}

\newcommand{\vp}{\varphi}
\newcommand{\lbal}{\left\{ \begin{array}{l}}
\newcommand{\lball}{\left\{ \begin{array}{ll}}
\newcommand{\ear}{\end{array} \right.}

\newcommand{\abs}{\\[5pt]}

\newcommand{\adb}{\allowdisplaybreaks}

\newcommand{\tme}{T_{max,\eps}}

\newcommand{\ueps}{u_\eps}
\newcommand{\veps}{v_\eps}

\newcommand{\heps}{h_\eps}
\newcommand{\feps}{f_\eps}

\newcommand{\yeps}{y_\eps}

\newcommand{\Teps}{\Theta_\eps}
\newcommand{\vepsx}{v_{\eps x}}
\newcommand{\vepsxx}{v_{\eps xx}}
\newcommand{\vepsxxx}{v_{\eps xxx}}
\newcommand{\vepsxxxx}{v_{\eps xxxx}}
\newcommand{\vepst}{v_{\eps t}}
\newcommand{\uepsx}{u_{\eps x}}
\newcommand{\uepsxx}{u_{\eps xx}}
\newcommand{\uepsxxx}{u_{\eps xxx}}
\newcommand{\uepst}{u_{\eps t}}
\newcommand{\Tepsx}{\Theta_{\eps x}}
\newcommand{\Tepsxx}{\Theta_{\eps xx}}
\newcommand{\Tepst}{\Theta_{\eps t}}
\newcommand{\gaeps}{\gamma_\eps}
\newcommand{\whv}{\wh{v}}
\newcommand{\whu}{\wh{u}}
\renewcommand{\div}{{\rm div} \,}
\newcommand{\nas}{\nabla^s}
\begin{document}
\adb
%
%
\title{
Large-data 
solutions in one-dimensional thermoviscoelasticity\\
involving temperature-dependent viscosities}
\author{
Michael Winkler\footnote{michael.winkler@math.uni-paderborn.de}\\
{\small Universit\"at Paderborn, Institut f\"ur Mathematik}\\
{\small 33098 Paderborn, Germany} }
\date{}
\maketitle
\begin{abstract}
\noindent 
An initial-boundary value problem for 
\bas
	\lball
	u_{tt} = \big(\gamma(\Theta) u_{xt}\big)_x + au_{xx} - \big(f(\Theta)\big)_x,
	\qquad & x\in\Om, \ t>0, \\[1mm]
	\Theta_t = \Theta_{xx} + \gamma(\Theta) u_{xt}^2 - f(\Theta) u_{xt},
	\qquad & x\in\Om, \ t>0, 
	\ear
\eas
is considered in an open bounded real interval $\Om$.
Under the assumption that $\gamma\in C^0([0,\infty))$ and $f\in C^0([0,\infty))$ are such that $f(0)=0$, and
$k_\gamma \le \gamma \le K_\gamma$ as well as
\bas
	|f(\xi)| \le K_f \cdot (\xi+1)^\al
	\qquad \mbox{for all } \xi\ge 0
\eas
with some $k_\gamma>0, K_\gamma>0, K_f>0$ and $\alpha<\frac{3}{2}$, for all suitably regular initial data 
of arbitrary size a statement on global existence of a global weak solution is derived.\abs
\end{abstract}
%
%
%
%
%
%
%
\newpage
\section{Introduction}\label{intro}
Evolution systems of the form 
\be{00}
	\lbal
	u_{tt} = \big(\gamma(\Theta,u_x) u_{xt}\big)_x + \big(a(\Theta,u_x) u_x)_x - \big(f(\Theta,u_x)\big)_x, \\[1mm]
	\Theta_t = \Theta_{xx} + \gamma(\Theta,u_x) u_{xt}^2 - f(\Theta,u_x) u_{xt},
	\ear
\ee
arise in the dynamical theory of thermoviscoelasticity in two relevant places.
Firstly, (\ref{00}) represents the classical model for viscosity-driven conversion of mechanical energy into heat
within solid materials (\cite{roubicek}, \cite{racke_zheng}). 
Secondly, can be viewed as a simplified description of thermoviscoelastic processes in contexts of piezoelectric materials,
where the generally vector-valued mechanical displacement $u$ and the temperature $\Theta$ are coupled to a divergence-free electric displacement field 
\be{01}
	D=\eps E+ e \nas u, 
\ee
with $\eps, e, E$ and $\nas u=\frac{1}{2} (\na u + (\na u)^T)$ denoting the 
permittivity matrix, the piezoelectric coupling tensor, the electric field strength and the usual symmetric gradient, respectively,
through an additional appearance of the last summand in the identity
\be{02}
	\rho u_{tt} = \div \big( d S_t + C S - \Theta CB \big) - \div \big(e^T E \big),
\ee
where $\rho$, $d$, $C$ and $B$ represent the material density, the viscosity tensor, the tensor of elastic parameters, 
and the thermal-dilation tensor.
In one-dimensional settings, namely, assuming constancy of the ingredients $\eps$, $e$ and $\rho$ allows for replacing $E$
in (\ref{02}) on the basis of (\ref{01}) and the solenoidality property $\div D=0$ to indeed obtain the first equation in (\ref{00})
with
\be{03}
	\gamma=\frac{d}{\rho},
	\qquad 
	a = \frac{C}{\rho} + \frac{e^2}{\eps\rho}
	\qquad \mbox{and} \qquad
	f=\frac{\Theta CB}{\rho}
\ee
(cf.~also \cite{fricke}).\abs
Subsequent to the seminal papers \cite{dafermos} and \cite{dafermos_hsiao_smooth}, 
the mathematical theory for the one-dimensional problem (\ref{00}) 
and close relatives has achieved a meanwhile rather well-developed state 
in scenarios in which its core ingredients $\gamma$ and $a$ do not explicitly depend on the temperature variable.
In \cite{racke_zheng}, for constant $\gamma$ and large classes of functions $a=a(u_x)$ and $F=F(u_x)$ 
in $f=\Theta F$ some weak solutions to an associated initial-boundary value problem have been shown to exist globally in time.
Results on global small-data solvability can be found in \cite{zheng_shen}, \cite{shibata} and \cite{kim},
and under appropriate though yet physically meaningful conditions also some globally smooth large-data
solutions have been constructed (\cite{guo_zhu}, \cite{shen_zheng_zhu}, \cite{chen_hoffmann}).
In the presence of such constant $\gamma$, even results on large-time stabilization toward stationary states are available
(\cite{racke_zheng}, \cite{hsiao_luo}).
Far-reaching results on global solvability extend to cases in which $\gamma$ is allowed to exhibit suitable dependencies on the
strain $u_x$ (\cite{watson}, \cite{jiang_QAM1993}).
Apart from that, a considerable literature has been concerned with higher-dimensional relatives of (\ref{00})
(see \cite{roubicek}, \cite{mielke_roubicek}, \cite{blanchard_guibe},
\cite{rossi_roubicek_interfaces13}, \cite{gawinecki_zajaczkowski_cpaa}, 
\cite{gawinecki_zajaczkowski}, \cite{roubicek_nodea2013},
\cite{pawlow_zajaczkowski_cpaa17} and \cite{owczarek_wielgos} for a small
selection), 
as well as with the viscosity-free classical models for thermoelasticity,
as reducing to (\ref{00}) with $\gamma\equiv 0$ in one-dimensional domains 
(cf.~e.g., \cite{slemrod}, \cite{bies_cieslak}, \cite{cieslak}, \cite{jiang1990}, \cite{racke90}, \cite{racke_shibata} 
and \cite{racke_shibata_zheng}).\abs
For versions of (\ref{00}) involving temperature dependencies of the core constituents $\gamma$ and $a$, however,
the knowledge seems significantly sparser, 
although recent experimental observations have hinted on some moderate but
non-negligible variations of $\gamma$ and $a$ in (\ref{00}) with respect to $\Theta$
(\cite{friesen}; see also \cite{Gubinyi2007}).
In fact, results on boundary value problems for (\ref{00}) with $\gamma=\gamma(\Theta), a=a(\Theta)$ and
$f=f(\Theta)$ apparently reduce to findings on local existence for initial data 
of arbitrary size, while solutions existing on large time intervals or even globally 
seem to have been constructed only in settings of suitably small data
(\cite{fricke}, \cite{meyer}, \cite{claes_win}, \cite{win_AMOP}
; see also \cite{claes_lankeit_win} for a result on global
generalized solvability in higher-dimensional cases when $f\equiv 0$).\abs
{\bf Main results.} \quad
The purpose of the present manuscript is to indicate that a considerably farther-reaching theory for (\ref{00}) 
with temperature-dependent ingredients can developed at least in cases when $a$ is constant.
Beyond the intention to address the corresponding exclusively thermo-mechanical scenarios from standard viscoelasticity
with $\gamma=\gamma(\Theta)$ and $a=const.$, this is motivated by
the ambition to achieve some basic understanding of the dynamics 
generated by piezoelectric interaction mechanisms such as those described near (\ref{01})-(\ref{03}).
Particularly in situations when, as suggested e.g.~by \cite{GutierrezLemini2014}, \cite{friesen} and \cite{Gubinyi2007}, 
$d=d(\Theta)$ and $C=C(\Theta)$ are linked via the relationship $d(\Theta)=\tau C(\Theta)$ 
with $\tau>0$ representing a certain retardation time constant quantifying mechanical losses, namely,
exemplary limit scenarios in (\ref{01}), (\ref{02}) and (\ref{03}) in which $\tau\to +\infty$ and $C\to 0$, but still $\tau C=O(1)$, propose to examine (\ref{00}) under the assumptions that
\be{04}
	\gamma\equiv \gamma(\Theta) = \frac{d(\Theta)}{\rho}
	\qquad \mbox{and} \qquad
	a = \frac{e^2}{\eps\rho} \equiv const.
\ee
Additionally assuming a general dependence of $f$ on $\Theta$, and supplementing (\ref{00}) by meaningful
boundary conditions, we are thus led to considering the initial-boundary value problem
\be{0}
	\lball
	u_{tt} = \big(\gamma(\Theta) u_{xt}\big)_x + au_{xx} - \big(f(\Theta)\big)_x,
	\qquad & x\in\Om, \ t>0, \\[1mm]
	\Theta_t = \Theta_{xx} + \gamma(\Theta) u_{xt}^2 - f(\Theta) u_{xt},
	\qquad & x\in\Om, \ t>0, \\[1mm]
	u=0, \quad \Theta_x=0,
	\qquad & x\in\pO, \ t>0, \\[1mm]
	u(x,0)=u_0(x), \quad u_t(x,0)=u_{0t}(x), \quad \Theta(x,0)=\Theta_0(x),
	\qquad & x\in\Om,
	\ear
\ee
with prescribed initial data $u_0, u_{0t}$ and $\Theta_0$ in an open bounded real interval.\abs
This problem will be addressed in the fairly natural framework of weak solvability specified as follows.
\begin{defi}\label{dw}
  Let $\Om\subset\R$ be a bounded open interval, let $\gamma\in C^0([0,\infty))$, $f\in C^0([0,\infty))$ and $a\in\R$, and let
  $u_0\in L^1(\Om)$, $u_{0t}\in L^1(\Om)$ and $\Theta_0\in L^1(\Om)$.
  Then by a {\em global weak solution} of (\ref{0}) we mean a pair $(u,\Theta)$ of functions
  \be{w1}
	\lbal
	u\in C^0([0,\infty);L^1(\Om)) \cap L^1_{loc}([0,\infty);W_0^{1,1}(\Om))
	\qquad \mbox{and} \\[1mm]
	\Theta\in L^1_{loc}([0,\infty);W^{1,1}(\Om))
	\ear
  \ee
  such that
  \be{w2}
	u_t \in L^2_{loc}(\Om\times [0,\infty))
  \ee
  as well as
  \be{w3}
	\big\{ \gamma(\Theta) u_{xt} \, , \, \gamma(\Theta) u_{xt}^2 \ , \, f(\Theta) \, , \, f(\Theta) u_{xt} \big\}
	\subset L^1_{loc}(\bom\times [0,\infty)),
  \ee
  that $u(\cdot,0)=u_0$ a.e.~in $\Om$ and $\Theta\ge 0$ a.e.~in $\Om\times (0,\infty)$, and that
  \bea{wu}
	- \int_0^\infty \io u_t \vp_t
	- \io u_{0t} \vp(\cdot,0)
	= - \int_0^\infty \io \gamma(\Theta) u_{xt} \vp_x
	- a \int_0^\infty \io u_x \vp_x
	+ \int_0^\infty \io f(\Theta) \vp_x
  \eea
  for all $\vp\in C_0^\infty(\Om\times [0,\infty))$, as well as
  \bea{wt}
	- \int_0^\infty \io \Theta \vp_t - \io \Theta_0 \vp(\cdot,0)
	= - \int_0^\infty \io \Theta_x \vp_x
	+ \int_0^\infty \io \gamma(\Theta) u_{xt}^2 \vp
	- \int_0^\infty \io f(\Theta) u_{xt} \vp
  \eea
  for each $\vp\in C_0^\infty(\bom\times [0,\infty))$.
\end{defi}
In this setting, our main results now assert global solvability under a suitable assumption on uniform non-degeneracy of 
$\gamma$ that seems consistent with the observations in \cite{friesen} and \cite{Gubinyi2007}, under a condition
on $f$ which is mild enough so as to include the standard linear dependence of thermal dilation mechanisms on temperature
(cf., e.g., \cite{roubicek}), and throughout a large class of initial data:
\begin{theo}\label{theo17}
  Let $\Om\subset\R$ be an open bounded interval, let $a>0$, and suppose that $\gamma\in C^0([0,\infty))$ and $f\in C^0([0,\infty))$
  are such that $f(0)=0$ and
  \be{g}
	k_\gamma \le \gamma(\xi) \le K_\gamma
	\qquad \mbox{for all } \xi\ge 0
  \ee
  and
  \be{f}
	|f(\xi)| \le K_f \cdot (\xi+1)^\al
	\qquad \mbox{for all } \xi\ge 0
  \ee
  with some $k_\gamma>0, K_\gamma>0, K_f>0$ and 
  \bas
	\al<\frac{3}{2}.
  \eas
  Then whenever
  \be{init}
	u_0\in W_0^{1,2}(\Om),
	\quad 
	u_{0t} \in L^2(\Om)
	\quad \mbox{and} \quad
	\Theta_0\in L^1(\Om)
  \ee
  are such that $\Theta_0\ge 0$ a.e.~in $\Om$, one can find functions
  \be{17.1}
	\lbal
	u \in C^0(\bom\times [0,\infty)) \cap L^\infty((0,\infty);W_0^{1,2}(\Om)) 
	\qquad \mbox{and} \\[1mm]
	\Theta\in L^\infty((0,\infty);L^1(\Om)) \cap \bigcap_{q\in [1,3)} L^q_{loc}(\bom\times [0,\infty))
		\cap \bigcap_{r\in [1,\frac{3}{2})} L^r_{loc}([0,\infty);W^{1,r}(\Om))
	\ear
  \ee
  such that
  \be{17.2}
	u_t \in L^\infty((0,\infty);L^2(\Om)) \cap L^2_{loc}([0,\infty);W_0^{1,2}(\Om)),
  \ee
  that $\Theta\ge 0$ a.e.~in $\Om\times (0,\infty)$, and that $(u,\Theta)$ forms a global weak solution of (\ref{0})
  in the sense of Definition \ref{dw}.
\end{theo}
{\bf Main ideas.} \quad
A key challenge for our analysis will consist in making appropriate use of the fundamental conservation property
expressed in the identity
\be{en}
	\frac{d}{dt} \bigg\{
	\frac{1}{2} \io v^2
	+ \frac{a}{2} \io u_x^2
	+ \io \Theta \bigg\}
	= 0,
	\qquad v=u_t,
\ee
that is formally associated with (\ref{0}), and that apparently marks a core difference to counterparts involving
functions $a=a(\Theta)$, for which only small-data solutions so far have been found to exist globally 
(\cite{fricke}, \cite{claes_win}).
In Section \ref{sect2}, a variant of (\ref{en}) will be derived for smooth solutions $(\veps,\ueps,\Teps)$
to a parabolic regularization of (\ref{0}) (see (\ref{0eps})), 
and Section \ref{sect3} will draw some preliminary and essentially immediate consequences
which entail global extensibility of these approximate solutions.
Section \ref{sect4} thereafter utilizes (\ref{en}) in a more substantial manner in order to reveal some
further regularity properties independent of the regularization parameter $\eps\in (0,1)$, seen to be valid under the
assumptions on $\gamma$ and $f$ made in Theorem \ref{theo17}.
Lemma \ref{lem11} will turn the collection of corresponding estimates into a statement on the existence of a limit object
$(v,u,\Theta)$ that is obtained on taking $\eps=\eps_j \searrow 0$ along a suitable sequence $(\eps_j)_{j\in\N}$,
and that can already at this point be seen to satisfy the first sub-problem of (\ref{0}) in the intended sense.\abs
The remaining part of our analysis will be devoted to the associated limit passage in the crucial quadratic contribution
to the weak formulations of the respective equations for temperature evolution, that is, to the derivation of the relation
\be{05}
	\int_0^\infty \io \gaeps(\Teps) \vepsx^2 \vp
	\to \int_0^\infty \io \gamma(\Theta) v_x^2 \vp
	\qquad \mbox{as } \eps=\eps_j\searrow 0
\ee
for arbitrary test functions $\vp\in C_0^\infty(\bom\times [0,\infty))$.
In its essence, this will amount to turning a basic weak $L^2$-approximation property of the $v_{\eps_j x}$
as having been obtained by a simple compactness argument in Lemma \ref{lem11}, into a corresponding strong convergence
feature. \abs
In Section \ref{sect5}, this will be achieved on the basis of a suitably careful testing procedure applied to the weak
identity (\ref{wu}). We will thereby exclude an unfavorable drop in the behavior of $L^2$ norms in this weak convergence process
by rigorously confirming a temporally localized variant of (\ref{en}) that is expressed in an inequality of the form
\be{06}
	\int_0^\infty \io \zeta \gamma(\Theta) v_x^2
	\ge \frac{a}{2} \int_0^\infty \io \zeta' u_x^2
	+ \frac{1}{2} \io u_{0t}^2
	+ \frac{a}{2} \io u_{0x}^2
	+ \frac{1}{2} \int_0^\infty \io \zeta' v^2
	+ \int_0^\infty \io \zeta f(\Theta) v_x
\ee
for each nonincreasing $\zeta\in C_0^\infty([0,\infty))$ fulfilling $\zeta(0)=1$ (Lemma \ref{lem15}).
This will entail that
\be{07}
	\sqrt{\gaeps(\Teps)} \vepsx \to \sqrt{\gamma(\Theta)} v_x
	\quad \mbox{in } L^2_{loc}(\bom\times [0,\infty))
	\qquad \mbox{as } \eps=\eps_j\searrow 0
\ee
(Lemma \ref{lem16}) and thus, by implying (\ref{05}), complete our argument.
\mysection{Preliminaries. An energy-consistent parabolic regularization}\label{sect2}
Our basic approach toward the construction of solutions will be motivated by the observation that
upon the substitution $v:=u_t$, the problem (\ref{0}) formally becomes equivalent to the system
\be{0v}
	\lball
	v_t = \big( \gamma(\Theta) v_x\big)_x + au_{xx} - \big( f(\Theta)\big)_x,
	\qquad & x\in\Om, \ t>0, \\[1mm]
	u_t = v,
	\qquad & x\in\Om, \ t>0, \\[1mm]
	\Theta_t = \Theta_{xx} + \gamma(\Theta) v_x^2 - f(\Theta) v_x,
	\qquad & x\in\Om, \ t>0, \\[1mm]
	v=0, \quad u=0, \quad \Theta_x=0,
	\qquad & x\in\pO, \ t>0, \\[1mm]
	v(x,0)=u_{0t}(x), \quad u(x,0)=u_0(x), \quad \Theta(x,0)=\Theta_0(x),
	\qquad & x\in\Om,
	\ear
\ee
which can be viewed as the limit case of a family of semilinear parabolic flows.
To substantiate this, in the framework of the assumptions from Theorem \ref{theo17}
we let $(v_{0\eps})_{\eps\in (0,1)} \subset C_0^\infty(\Om)$,
$(u_{0\eps})_{\eps\in (0,1)} \subset C_0^\infty(\Om)$
and $(\Theta_{0\eps})_{\eps\in (0,1)} \subset C^\infty(\bom)$
be such that $\Theta_{0\eps} \ge 0$ in $\Om$ for all $\eps\in (0,1)$, and that as $\eps\searrow 0$ we have
\be{ie}
	\lbal
	v_{0\eps} \to u_{0t}
	\qquad \mbox{in } L^2(\Om), \\[1mm]
	u_{0\eps} \to u_0
	\qquad \mbox{in } W^{1,2}(\Om)
	\qquad \mbox{and} \\[1mm]
	\Theta_{0\eps} \to \Theta_0
	\qquad \mbox{in } L^1(\Om).
	\ear
\ee
We moreover let $(\gaeps)_{\eps\in (0,1)} \subset C^\infty([0,\infty))$ and 
$(\feps)_{\eps\in (0,1)} \subset C^\infty([0,\infty)) \cap L^\infty([0,\infty))$
be such that
\be{gaeps}
	k_\gamma \le \gaeps(\xi) \le K_\gamma(\xi)
	\qquad \mbox{for all $\xi\ge 0$ and } \eps\in (0,1),
\ee
that
\be{feps}
	|\feps(\xi)|\le K_f (\xi+1)^\al
	\qquad \mbox{for all $\xi\ge 0$ and } \eps\in (0,1),
\ee
and that
\be{gflimit}
	\gaeps\to \gamma
	\quad \mbox{and} \quad
	\feps\to f
	\quad \mbox{in } L^\infty_{loc}([0,\infty))
	\qquad \mbox{as } \eps\searrow 0,
\ee
and for $\eps\in (0,1)$ we consider
\be{0eps}
	\lball
	\vepst = - \eps \vepsxxxx + \big( \gaeps(\Teps) \vepsx \big)_x + a\uepsxx - \big( \feps(\Teps)\big)_x,
	\qquad & x\in\Om, \ t>0, \\[1mm]
	\uepst = \eps \uepsxx + \veps,
	\qquad & x\in\Om, \ t>0, \\[1mm]
	\Tepst = \Tepsxx + \gaeps(\Teps) \vepsx^2 - \feps(\Teps) \vepsx,
	\qquad & x\in\Om, \ t>0, \\[1mm]
	\veps=\vepsxx=0, \quad \ueps=0, \quad \Tepsx=0,
	\qquad & x\in\pO, \ t>0, \\[1mm]
	\veps(x,0)=v_{0\eps}(x), \quad \ueps(x,0)=u_{0\eps}(x), \quad \Teps(x,0)=\Theta_{0\eps}(x),
	\qquad & x\in\Om.
	\ear
\ee
Indeed, all these problems are accessible to standard theory of local solvability in parabolic systems:
\begin{lem}\label{lem1}
  Let $\eps\in (0,1)$.
  Then there exist $\tme\in (0,\infty]$ and functions
  \bas
	\lbal
	\veps\in C^0(\bom\times [0,\tme))
		\cap C^{4,1}(\bom\times (0,\tme)), \\[1mm]
	\ueps\in C^0(\bom\times [0,\tme))
		\cap C^{2,1}(\bom\times (0,\tme))
		\cap C^0([0,\tme);W_0^{1,2}(\Om)) \qquad \mbox{and} \\[1mm]
	\Teps \in C^0(\bom\times [0,\tme))
		\cap C^{2,1}(\bom\times (0,\tme))
	\ear
  \eas
  such that $\Teps\ge 0$ in $\bom\times [0,\tme)$, that (\ref{0eps}) is solved in the classical sense in 
  $\Om\times (0,\tme)$, and that
  \bea{ext}
	& & \hs{-15mm}
	\mbox{if $\tme<\infty$, \quad then } \nn\\
	& & \hs{-10mm}
	\limsup_{t\nearrow \tme} \Big\{
	\|\veps(\cdot,t)\|_{W^{2+2\eta,\infty}(\Om)}
	+ \|\ueps(\cdot,t)\|_{W^{1+\eta,\infty}(\Om)}
	+ \|\Teps(\cdot,t)\|_{W^{1+\eta,\infty}(\Om)}
	\Big\} = \infty
	\ \mbox{for all } \eta>0.
  \eea
\end{lem}
\proof
  This can be derived in a straightforward manner from the general results on existence and extensibility stated
  in \cite[Theorem 12.1, Theorem 12.5]{amann}; in a setting closely related to that in (\ref{0eps}), details on 
  a possible reduction to the framework from \cite{amann} can be found in \cite[Lemma 2.3]{claes_lankeit_win}.
\qed
Now three elementary testing procedures show that the regularization made in (\ref{0eps}) accords well with the 
energy structure in (\ref{en}).
The outcomes of the first two of these will be of independent use 
(see Lemma \ref{lem8} and Lemma \ref{lem16} below), and will thus be
recorded here separately:
\begin{lem}\label{lem02}
  If $\eps\in (0,1)$, then
  \be{02.1}
	\frac{1}{2} \frac{d}{dt} \io \veps^2 
	+ \frac{a}{2} \frac{d}{dt} \io \uepsx^2
	+ \io \gaeps(\Teps) \vepsx^2
	+ \eps \io \vepsxx^2
	+ \eps a \io \uepsxx^2
	= \io f(\Teps) \vepsx
	\qquad \mbox{for all } t\in (0,\tme).
  \ee
\end{lem}
\proof
  According to the first equation in (\ref{0eps}), integrations by parts show that
  \bas
	\frac{1}{2} \frac{d}{dt} \io \veps^2
	&=& - \eps \io \veps \vepsxxxx
	+ \io \veps \cdot \big( \gaeps(\Teps) \vepsx + a\uepsx - \feps(\Teps)\big)_x \\
	&=& - \eps \io \vepsxx^2
	- \io \gaeps(\Teps) \vepsx^2
	- a \io \uepsx \vepsx
	+ \io f(\Teps) \vepsx
	\qquad \mbox{for all } t\in (0,\tme).
  \eas
  Therefore, (\ref{02.1}) follows upon observing that by the second equation in (\ref{0eps}),
  \bas
	\frac{a}{2} \frac{d}{dt} \io \uepsx^2
	= a \io \uepsx \cdot \big( \eps\uepsxx + \veps\big)_x 
	= - \eps a \io \uepsxx^2
	+ a \io \uepsx \vepsx
  \eas
  for all $t\in (0,\tme)$.
\qed
When supplemented by a simple observation on evolution of the functional $\io \Teps$, the above yields
the following regularized counterpart of (\ref{en}).
\begin{lem}\label{lem2}
  The solution of (\ref{0eps}) satisfies
  \be{2.1}
	\frac{d}{dt} \bigg\{
	\frac{1}{2} \io \veps^2
	+ \frac{a}{2} \io \uepsx^2
	+ \io \Teps \bigg\}
	+ \eps \io \vepsxx^2
	+ \eps a \io \uepsxx^2 =0
	\qquad \mbox{for all } t\in (0,\tme).
  \ee
\end{lem}
\proof
  As
  \bas
	\frac{d}{dt} \io \Teps
	= \io \gaeps(\Teps) \vepsx^2
	- \io \feps(\Teps) \vepsx
	\qquad \mbox{for all } t\in (0,\tme)
  \eas
  by the third equation in (\ref{0eps}), we obtain (\ref{2.1}) from (\ref{02.1}) thanks to two favorable cancellations.
\qed
Evident consequences of (\ref{2.1}) provide some fundamental regularity information:
\begin{lem}\label{lem3}
  There exists $C>0$ such that
  \be{3.1}
	\io \veps^2(\cdot,t) \le C
	\qquad \mbox{for all } t\in (0,\tme)
  \ee
  and
  \be{3.2}
	\io \uepsx^2(\cdot,t) \le C
	\qquad \mbox{for all } t\in (0,\tme)
  \ee
  and
  \be{3.3}
	\io \Teps \le C
	\qquad \mbox{for all } t\in (0,\tme)
  \ee
  as well as
  \be{3.4}
	\eps \int_0^t \io \vepsxx^2 \le C
	\qquad \mbox{for all } t\in (0,\tme)
  \ee
  and
  \be{3.5}
	\eps \int_0^t \io \uepsxx^2 \le C
	\qquad \mbox{for all } t\in (0,\tme).
  \ee
\end{lem}
\proof
  An integration of (\ref{2.1}) results in the identity
  \bas
	& & \hs{-20mm}
	\frac{1}{2} \io \veps^2 
	+ \frac{a}{2} \io \uepsx^2
	+ \io \Teps
	+ \eps \int_0^t \io \vepsxx^2
	+ \eps a \int_0^t \io \uepsxx^2 \\
	&=& \frac{1}{2} \io v_{0\eps}^2 
	+ \frac{a}{2} \io u_{0\eps x}^2
	+ \io \Theta_{0\eps}
  \eas
  for all $t\in (0,\tme)$ and $\eps\in (0,1)$.
  This already implies (\ref{3.1})-(\ref{3.5}) due to the nonnegativity of $\Teps$ and the fact that
  $\sup_{\eps\in (0,1)} \io v_{0\eps}^2$,
  $\sup_{\eps\in (0,1)} \io u_{0\eps x}^2$
  and
  $\sup_{\eps\in (0,1)} \io \Theta_{0\eps}$
  are all finite due to (\ref{ie}).
\qed
\mysection{Global solvability of the approximate problems}\label{sect3}
The goal of this section is to make sure that for each fixed $\eps\in (0,1)$,
the first alternative in the extensibility criterion cannot occur.
Accordingly focusing on higher regularity properties possibly depending on $\eps$ here, in a first step in this regard
we may rely on the artificial fourth-order diffusion mechanism in (\ref{0eps}), and on our assumption that
$(\feps)_{\eps\in (0,1)} \subset L^\infty([0,\infty))$, to obtain the following.
\begin{lem}\label{lem31}
  If $\tme<\infty$ for some $\eps\in (0,1)$, then there exists $C(\eps)>0$ such that
  \be{31.1}
	\io \vepsx^2(\cdot,t) \le C(\eps)
	\qquad \mbox{for all } t\in (0,\tme)
  \ee
  and
  \be{31.2}
	\io \uepsxx^2(\cdot,t) \le C(\eps)
	\qquad \mbox{for all } t\in (0,\tme),
  \ee
  and moreover we have
  \be{31.3}
	\int_0^{\tme} \io \vepsxxx^2 < \infty.
  \ee
\end{lem}
\proof
  Relying on the smoothness properties of $(\veps,\ueps,\Teps)$ recorded in Lemma \ref{lem1}, we may use $-\vepsxx$ as a test
  function in the first equation from (\ref{0eps}) to see that
  for all $t\in (0,\tme)$,
  \bea{31.4}
	\frac{1}{2} \frac{d}{dt} \io \vepsx^2
	+ \eps \io \vepsxxx^2
	&=& \io \gaeps(\Teps) \vepsx \vepsxxx
	- a \io \uepsxx \vepsxx
	- \io \feps(\Teps) \vepsxxx,
  \eea
  because $\vepsxx=0$ on $\pO\times (0,\tme)$.
  Here, by Young's inequality and (\ref{gaeps}),
  \bas
	\io \gaeps(\Teps)\vepsx\vepsxxx
	\le \frac{\eps}{4} \io \vepsxxx^2
	+ \frac{K_\gamma^2}{\eps} \io \vepsx^2
	\qquad \mbox{for all } t\in (0,\tme),
  \eas
  while using that $c_1\equiv c_1(\eps):=\|\feps\|_{L^\infty((0,\infty))}$ is finite, by means of Young's inequality we moreover
  find that
  \bas
	- \io \feps(\Teps) \vepsxxx
	\le \frac{\eps}{4} \io \vepsxxx^2
	+ \frac{c_1^2 |\Om|}{\eps}
	\qquad \mbox{for all } t\in (0,\tme).
  \eas
  As furthermore
  \bas
	- a \io \uepsxx \vepsxx
	&=& - a \io \uepsxx \cdot \big\{ u_{\eps xxt} - \eps u_{\eps xxxx} \big\} \\
	&=& - \frac{a}{2} \frac{d}{dt} \io \uepsxx^2
	- \eps a \io \uepsxxx^2 \\
	&\le& - \frac{a}{2} \frac{d}{dt} \io \uepsxx^2
	\qquad \mbox{for all } t\in (0,\tme)
  \eas
  by the second equation in (\ref{0eps}) and the fact that $\uepsxx=\frac{1}{\eps} (\uepst-\veps)=0$ on $\pO\times (0,\tme)$,
  from (\ref{31.4}) we infer that $\yeps(t):=1+\io \vepsx^2(\cdot,t) + a \io \uepsxx^2(\cdot,t)$, $t\in [0,\tme)$,
  satisfies
  \be{31.5}
	\yeps'(t) + \eps \io \vepsxxx^2 \le c_2 \yeps(t)
	\qquad \mbox{for all } t\in (0,\tme)
  \ee
  with $c_2\equiv c_2(\eps):=\max\big\{ \frac{2K_\gamma^2}{\eps} \, , \, \frac{2c_1^2 |\Om|}{\eps}\big\}$.
  By an ODE comparison argument, this firstly implies that $\yeps(t)\le c_3\equiv c_3(\eps):=\yeps(0) e^{c_2\tme}$
  for all $t\in (0,\tme)$, and hence entails both (\ref{31.1}) and (\ref{31.2}), 
  and thereupon an integration in (\ref{31.5}) shows that $\eps\int_0^{\tme} \io \vepsxxx^2 \le \yeps(0) + c_2 c_3 \tme$,
  and that thus also (\ref{31.3}) holds.
\qed
Again thanks to the boundedness of each individual $\feps$, by means of heat semigroup estimates the above can be turned
into a time-independent estimate for $\Teps$ within a topological framework already compatible with that appearing in 
(\ref{ext}).
\begin{lem}\label{lem32}
  Suppose that $\tme<\infty$ for some $\eps\in (0,1)$.
  Then given any $\sig\in (1,\frac{5}{4})$, one can find $C(\sig,\eps)>0$ such that
  \be{32.1}
	\|\Teps(\cdot,t)\|_{W^{\sig,\infty}(\Om)} \le C(\sig,\eps)
	\qquad \mbox{for all } t\in (0,\tme).
  \ee
\end{lem}
\proof
  Since $\sig<\frac{5}{4}$, we can choose $\beta\in (\frac{3}{4},\frac{7}{8})$ suitably close to $\frac{7}{8}$ such that
  $2\beta-\frac{1}{2}>\sig$, and letting $A$ denote the self-adjoint
  realization of $-(\cdot)_{xx} +1$ under homogeneous Neumann boundary conditions in $L^2(\Om)$,
  with domain given by $D(A):=\{\psi\in W^{2,2}(\Om) \ | \ \frac{\pa\psi}{\pa\nu}=0 \mbox{ on } \pO\}$,
  we note that its fractional power $A^\beta$ satisfies $D(A^\beta) \hra W^{\sig,\infty}(\Om)$ (\cite{henry}),
  whence there exists $c_1>0$ such that
  \be{32.02}
	\|\psi\|_{W^{\sig,\infty}(\Om)} + \|\psi\|_{L^\infty(\Om)} \le c_1 \|A^\beta \psi\|_{L^2(\Om)}
	\qquad \mbox{for all } \psi\in D(A^\beta).
  \ee
  Now the basis of a Duhamel representation associated with the identity
  \bas
	\Tepst + A\Teps = \heps := \gaeps(\Teps) \vepsx^2 - \feps(\Teps) \vepsx + \Teps,
  \eas
  using known smoothing properties of the corresponding analytical semigroup $(e^{-tA})_{t\ge 0}$ we can find $c_2>0$ such that
  \bea{32.2}
	\hs{-8mm}
	\|A^\beta \Teps(\cdot,t)\|_{L^2(\Om)}
	&=& \bigg\| e^{-tA} A^\beta \Theta_{0\eps} + \int_0^t A^\beta e^{-(t-s)A} \heps(\cdot,s) ds \bigg\|_{L^2(\Om)} \nn\\
	&\le& c_2 \|A^\beta \Theta_{0\eps}\|_{L^2(\Om)} 
	+ c_2 \int_0^t (t-s)^{-\beta} \|\heps(\cdot,s)\|_{L^2(\Om)} ds
	\quad \mbox{for all } t\in (0,\tme).
  \eea
  Here, according to the boundedness of $\gaeps$ and $\feps$ on $[0,\infty)$ we can pick $c_3=c_3(\eps)>0$ such that
  \bas
	|\heps| \le c_3 \vepsx^2 + \Teps + c_3
	\qquad \mbox{in } \Om\times (0,\tme),
  \eas
  and observing that by a Gagliardo-Nirenberg inequality and (\ref{31.1}) there exist $c_4=c_4(\eps)>0$ and $c_5=c_5(\eps)>0$
  fulfilling
  \bas
	\|\vepsx^2\|_{L^2(\Om)}
	&=& \|\vepsx\|_{L^4(\Om)}^2
	\le c_4\|\vepsxxx\|_{L^2(\Om)}^\frac{1}{4} \|\vepsx\|_{L^2(\Om)}^\frac{7}{4} + c_4 \|\vepsx\|_{L^2(\Om)}^2 \\
	&\le& c_5\|\vepsxxx\|_{L^2(\Om)}^\frac{1}{4} + c_5
	\qquad \mbox{for all } t\in (0,\tme),
  \eas
  and that due to the H\"older inequality, (\ref{3.3}) and (\ref{32.02}) we have
  \bas
	\|\Teps\|_{L^2(\Om)} \le \|\Teps\|_{L^\infty(\Om)}^\frac{1}{2} \|\Teps\|_{L^1(\Om)}^\frac{1}{2}
	\le c_6 \|\Teps\|_{L^\infty(\Om)}^\frac{1}{2}
	\le c_1^\frac{1}{2} c_6 \|A^\beta \Teps\|_{L^2(\Om)}^\frac{1}{2}
	\qquad \mbox{for all } t\in (0,\tme)
  \eas
  with some $c_6>0$, we can rely on Young's inequality and (\ref{32.02}) to estimate
  \bea{32.3}
	c_2 \int_0^t (t-s)^{-\beta} \|\heps(\cdot,s)\|_{L^2(\Om)} ds
	&\le& c_2 c_3 \int_0^t (t-s)^{-\beta} \|\vepsx^2(\cdot,s)\|_{L^2(\Om)} ds \nn\\
	& & + c_2 \int_0^t (t-s)^{-\beta} \|\Teps(\cdot,s)\|_{L^2(\Om)} ds 
	+ c_2 c_3 \int_0^t (t-s)^{-\beta} ds \nn\\
	&\le& c_2 c_3 c_5 \int_0^t (t-s)^{-\beta} \cdot \Big\{ \|\vepsxxx(\cdot,s)\|_{L^2(\Om)}^\frac{1}{4} + 1 \Big\} ds \nn\\
	& & + c_1^\frac{1}{2} c_2 c_6 \int_0^t (t-s)^{-\beta} \|A^\beta \Teps(\cdot,s)\|_{L^2(\Om)}^\frac{1}{2} ds
	+ c_2 c_3 \int_0^t (t-s)^{-\beta} ds \nn\\
	&\le& c_2 c_3 c_5 \int_0^t \io \vepsxxx^2
	+ c_2 c_3 c_5 \int_0^t (t-s)^{-\frac{8\beta}{7}} ds \nn\\
	& & + c_1^\frac{1}{2} c_2 c_6 M_\eps^\frac{1}{2}(T) \cdot \int_0^t (t-s)^{-\beta} ds 
	+ (c_2 c_3 c_5 + c_2 c_3) \int_0^t (t-s)^{-\beta} ds \nn\\[2mm]
	& & \hs{26mm}
	\qquad \mbox{for all $t\in (0,T)$ and } T\in (0,\tme),
  \eea
  where we have set
  \be{32.4}
	M_\eps(T):=\sup_{t\in (0,T)} \|A^\beta \Teps(\cdot,t)\|_{L^2(\Om)}
	\qquad \mbox{for } T\in (0,\tme).
  \ee
  As $\beta<\frac{8\beta}{7}<1$, a combination of (\ref{32.3}), (\ref{31.3}) and (\ref{32.2}) as well as the fact that 
  $\Theta_{0\eps} \in C_0^\infty(\Om) \subset D(A^\beta)$ thus reveals the existence of $c_7=c_7(\eps)>0$ such that
  \bas
	\|A^\beta \Teps(\cdot,t)\|_{L^2(\Om)}
	\le c_7 M_\eps^\frac{1}{2}(T) + c_7
	\qquad \mbox{for all $t\in (0,T)$ and } T\in (0,\tme),
  \eas
  and that therefore, by Young's inequality,
  \bas
	M_\eps(T) \le c_7 M_\eps^\frac{1}{2}(T) + c_7
	\le \frac{1}{2} M_\eps(T) + \frac{c_7^2}{2} + c_7
	\qquad \mbox{for all } T\in (0,\tme).
  \eas
  Hence, $M_\eps(T)\le c_7^2 + 2c_7$ for all $T\in (0,\tme)$, so that (\ref{32.1}) results from (\ref{32.4}) and (\ref{32.02}).
\qed
Along with (\ref{31.2}), this in turn implies regularity properties of the three rightmost summands in the 
first equation in (\ref{0eps}) which are sufficient to ensure bounds also for $\veps$ with respect to some of the norms
appearing in (\ref{0eps}).
\begin{lem}\label{lem33}
  If $\eps\in (0,1)$ is such that $\tme<\infty$, then for each $\sig\in (3,\frac{7}{2})$ there exists $C(\sig,\eps)>0$ fulfilling
  \be{33.1}
	\|\veps(\cdot,t)\|_{W^{\sig,\infty}(\Om)} \le C(\sig,\eps)
	\qquad \mbox{for all } t\in (0,\tme).
  \ee
\end{lem}
\proof
  Given $\sig\in (3,\frac{7}{2})$, we can fix $\beta\in (\frac{3}{4},1)$ in such a way that $4\beta-\frac{1}{2}>\sig$,
  meaning that for the self-adjoint operator $A:=\eps(\cdot)_{xxxx}$ in $L^2(\Om)$ with domain 
  $D(A):=\{\psi\in W^{4,2}(\Om) \ | \ \psi=\psi_{xx}=0 \mbox{ on } \pO\}$, the corresponding fractional power 
  $A^\beta$ satisfies $D(A^\beta) \hra W^{\sig,\infty}(\Om)$ (\cite{henry}), whence with some $c_1=c_1(\eps)>0$ we have
  \be{33.2}
	\|\psi\|_{W^{\sig,\infty}(\Om)} + \|\psi_{xxx}\|_{L^\infty(\Om)} \le c_1 \|A^\beta \psi\|_{L^2(\Om)}
	\qquad \mbox{for all } \psi\in D(A^\beta).
  \ee
  For $T\in (0,\tme)$, we now estimate
  \be{33.03}
	M_\eps(T):=\sup_{t\in (0,T)} \|A^\beta \veps(\cdot,t)\|_{L^2(\Om)}
  \ee
  by observing that in the identity
  \be{33.3}
	\vepst + A\veps = \heps:=\gaeps(\Teps) \vepsxx + \gaeps'(\Teps) \Tepsx \vepsx + a\uepsxx - \feps'(\Teps) \Tepsx,
  \ee
  Lemma \ref{lem32} and the continuity of $\gaeps,\gaeps'$ and $\feps'$ on $[0,\infty)$ ensure the existence of $c_2=c_2(\eps)>0$
  satisfying
  \bas
	\|\heps\|_{L^2(\Om)}
	\le c_2\|\vepsxx\|_{L^2(\Om)}
	+ c_2\|\vepsx\|_{L^2(\Om)}
	+ c_2 \|\uepsxx\|_{L^2(\Om)}
	+ c_2
	\qquad \mbox{for all } t\in (0,\tme).
  \eas
  Since the Gagliardo-Nirenberg inequality provides $c_3>0$ such that
  \bas
	\|\vepsxx\|_{L^2(\Om)}
	\le c_3\|\vepsxxx\|_{L^2(\Om)}^\frac{1}{2} \|\vepsx\|_{L^2(\Om)}^\frac{1}{2} + c_3\|\vepsx\|_{L^2(\Om)}
	\qquad \mbox{for all } t\in (0,\tme),
  \eas
  in view of (\ref{31.1}), (\ref{31.2}) and (\ref{33.2}) we thereby see that with some $c_4=c_4(\eps)>0$,
  \bas
	\|\heps\|_{L^2(\Om)}
	\le c_4\|\vepsxxx\|_{L^2(\Om)}^\frac{1}{2} + c_4
	\le c_1^\frac{1}{2} c_4 \|A^\beta \veps\|_{L^2(\Om)}^\frac{1}{2} + c_4
	\qquad \mbox{for all } t\in (0,\tme).
  \eas
  By means of a standard smoothing estimate for the semigroup $(e^{-tA})_{t\ge 0}$, we thus obtain that with some 
  $c_5=c_5(\eps)>0$,
  \bas
	\|A^\beta \veps(\cdot,t)\|_{L^2(\Om)}
	&=& \bigg\| e^{-tA} A^\beta v_{0\eps} + \int_0^t A^\beta e^{-(t-s)A} \heps(\cdot,s) ds \bigg\|_{L^2(\Om)} \\
	&\le& c_5 \|A^\beta v_{0\eps}\|_{L^2(\Om)}
	+ c_5 \int_0^t (t-s)^{-\beta} \|\heps(\cdot,s)\|_{L^2(\Om)} ds \\
	&\le& c_5 \|A^\beta v_{0\eps}\|_{L^2(\Om)}
	+ c_1^\frac{1}{2} c_4 c_5 M_\eps^\frac{1}{2}(T) \cdot \int_0^t (t-s)^{-\beta} ds \nn\\
	& & + c_4 c_5 \int_0^t (t-s)^{-\beta} ds
	\qquad \mbox{for all $t\in (0,T)$ and } T\in (0,\tme).
  \eas
  Since $\beta<1$ and $v_{0\eps}\in D(A^\beta)$, this yields $c_6=c_6(\eps)>0$ such that
  \bas
	\|A^\beta \veps(\cdot,t)\|_{L^2(\Om)}
	\le c_6 M_\eps^\frac{1}{2}(T) + c_6
	\qquad \mbox{for all $t\in (0,T)$ and } T\in (0,\tme),
  \eas
  so that the claim can readily be deduced on taking suprema here, and recalling (\ref{33.03}) and (\ref{33.2}).
\qed
In conclusion, each of our approximate solutions can actually be extended so as to exist globally:
\begin{lem}\label{lem4}
  For each $\eps\in (0,1)$, the solution of (\ref{0eps}) is global in time; that is, in Lemma \ref{lem1}
  we have $\tme=\infty$.
\end{lem}
\proof
  As $W^{2,2}(\Om) \hra W^{\sig,\infty}(\Om)$ whenever $\sig\in (1,\frac{3}{2})$ (\cite{henry}), 
  Lemma \ref{lem33} in conjunction with (\ref{31.2}) and Lemma \ref{lem32} shows that if $\tme$ was finite for some
  $\eps\in (0,1)$, then (\ref{ext}) would be violated. The claim thus results from Lemma \ref{lem1}.
\qed

\mysection{Further $\eps$-independent estimates. Construction of a limit triple}\label{sect4}	
Next addressing a key issue related to the construction of limit objects $v, u$ and $\Theta$ through an appropriate
extraction of subsequences of solutions to (\ref{0eps}), 
in this section we return to the consequences of the energy structure from (\ref{en}) formulated in Lemma \ref{lem3},
and plan to derive further $\eps$-independent regularity properties from these.\abs
We first prepare an exploitation of the heat diffusion mechanism in (\ref{0eps}) by stating the outcome
of a Gagliardo-Nirenberg interpolation, to be used both in Lemma \ref{lem5} and in Lemma \ref{lem6}:
\begin{lem}\label{lem44}
  Let $p>0$ and $q\in\R$. Then there exists $C=C(p,q)>0$ such that
  \be{44.1}
	\io (\Teps+1)^q
	\le C \cdot \bigg\{ \io (\Teps+1)^{p-2} \Tepsx^2 \bigg\}^\frac{q-1}{p+1}
	+ C
	\qquad \mbox{for all $t>0$ and } \eps\in (0,1).
  \ee
\end{lem}
\proof
  If $q>1$, then we employ the Gagliardo-Nirenberg inequality to find $c_1=c_1(p,q)>0$ fulfilling
  \bas
	\|\psi\|_{L^\frac{2q}{p}(\Om)}^\frac{2q}{p}
	\le c_1\|\psi_x\|_{L^2(\Om)}^\frac{2(q-1)}{p+1} \|\psi\|_{L^\frac{2}{p}(\Om)}^\frac{2(p+q)}{p(p+1)}
	+ c_1\|\psi\|_{L^\frac{2}{p}(\Om)}^\frac{2q}{p}
	\qquad \mbox{for all } \psi\in W^{1,2}(\Om).
  \eas
  An application to $\psi:=(\Teps+1)^\frac{p}{2}$ for $t>0$ and $\eps\in (0,1)$ then yields 
  \bas
	\io (\Teps+1)^q
	&\le& c_1\cdot \bigg\{ \io \Big| \big( (\Teps+1)^\frac{p}{2}\big)_x \Big|^2 \bigg\}^\frac{q-1}{p+1} 
		\cdot \bigg\{ \io (\Teps+1) \bigg\}^\frac{p+q}{p+1} \\
	& & + c_1 \cdot \bigg\{ \io (\Teps+1) \bigg\}^q
	\qquad \mbox{for all $t>0$ and } \eps\in (0,1),
  \eas
  whence in view of the evident identities
  \bas
	\io \Big| \big( (\Teps+1)^\frac{p}{2}\big)_x \Big|^2
	= \frac{p^2}{4} \io (\Teps+1)^{p-2} \Tepsx^2
  \eas
  and $\io (\Teps+1)= |\Om| + \io \Teps$, the claimed inequality in this case results from (\ref{3.3}).\abs
  If $q\le 1$, however, then simply using that $(\xi+1)^q \le \xi+2$ for all $\xi\ge 0$ in estimating
  \bas
	\io (\Teps+1)^q 
	\le \io (\Teps+2)
	\le 2|\Om| + \io \Teps
	\qquad \mbox{for all $t>0$ and } \eps\in (0,1)
  \eas
  we obtain (\ref{44.1}) even as an immediate consequence of (\ref{3.3}).
\qed
Now thanks to our overall assumption that $f$ satisfy (\ref{f}) with $\al<\frac{3}{2}$,
this interpolation result can be used to control ill-signed contributions to the evolution of the 
functionals $\io (\Teps+1)^p$ in their concave range where $p\in (0,1)$. 
We thereby obtain the following boundedness property of the quantities accordingly dissipated due to diffusion:
\begin{lem}\label{lem5}
  For each $p\in (0,1)$ and any $T>0$ there exists $C(p,T)>0$ with the property that
  \be{5.1}
	\int_0^T \io (\Teps+1)^{p-2} \Tepsx^2
	\le C(p,T)
	\qquad \mbox{for all } \eps\in (0,1).
  \ee
\end{lem}
\proof
  We once again use the third equation in (\ref{0eps}) along with the homogeneous Neumann boundary conditions for $\Teps$
  to see on the basis of an integration by parts that
  \bea{5.2}
	- \frac{1}{p} \frac{d}{dt} \io (\Teps+1)^p
	&=& - \io (\Teps+1)^{p-1} \cdot \big\{ \Tepsxx + \gaeps(\Teps) \vepsx^2 - \feps(\Teps) \vepsx \big\} \nn\\
	&=& -(1-p) \io (\Teps+1)^{p-2} \Tepsx^2
	- \io (\Teps+1)^{p-1} \gaeps(\Teps) \vepsx^2 \nn\\
	& & + \io (\Teps+1)^{p-1} f(\Teps) \vepsx
	\qquad \mbox{for all } t>0,
  \eea
  where using (\ref{gaeps}) we can estimate
  \be{5.3}
	\io (\Teps+1)^{p-1} \gaeps(\Teps) \vepsx^2
	\ge k_\gamma \io (\Teps+1)^{p-1} \vepsx^2
	\qquad \mbox{for all } t>0,
  \ee
  while by Young's inequality and (\ref{feps}),
  \bea{5.4}
	\io (\Teps+1)^{p-1} \feps(\Teps) \vepsx
	&\le& k_\gamma \io (\Teps+1)^{p-1} \vepsx^2
	+ \frac{1}{4 k_\gamma} \io (\Teps+1)^{p-1} \feps^2(\Teps) \nn\\
	&\le& k_\gamma \io (\Teps+1)^{p-1} \vepsx^2
	+ \frac{K_f^2}{4 k_\gamma} \io (\Teps+1)^{p+2\al-1}
	\qquad \mbox{for all } t>0.
  \eea
  Now Lemma \ref{lem44} provides $c_1=c_1(p)>0$ such that
  \be{5.5}
	\frac{K_f^2}{4 k_\gamma} \io (\Teps+1)^{p+2\al-1}
	\le c_1 \cdot \bigg\{ \io (\Teps+1)^{p-2} \Tepsx^2 \bigg\}^\frac{p+2\al-2}{p+1} + c_1
	\qquad \mbox{for all $t>0$ and } \eps\in (0,1),
  \ee
  where our overall assumption that $\al<\frac{3}{2}$ asserts that $\frac{p+2\al-2}{p+1}<1$.
  As $p<1$, we may therefore invoke Young's inequality once more to find $c_2=c_2(p)>0$ fulfilling
  \bas
	c_1 \cdot \bigg\{ \io (\Teps+1)^{p-2} \Tepsx^2 \bigg\}^\frac{p+2\al-2}{p+1} 
	\le \frac{1-p}{2} \io (\Teps+1)^{p-2} \Tepsx^2
	+ c_2
	\qquad \mbox{for all $t>0$ and } \eps\in (0,1).
  \eas
  When combined with (\ref{5.5}), (\ref{5.4}), (\ref{5.3}) and (\ref{5.2}), this shows that
  \bas
	- \frac{1}{p} \frac{d}{dt} \io (\Teps+1)^p
	+ \frac{1-p}{2} \io (\Teps+1)^{p-2} \Tepsx^2
	\le c_1 + c_2
	\qquad \mbox{for all $t>0$ and } \eps\in (0,1),
  \eas
  and that thus
  \be{5.6}
	\frac{1-p}{2} \int_0^T \io (\Teps+1)^{p-2} \Tepsx^2
	\le \frac{1}{p} \io \big(\Teps(\cdot,T)+1\big)^p + (c_1+c_2)T
	\qquad \mbox{for all $T>0$ and } \eps\in (0,1).
  \ee
  Again using that $p<1$ in estimating
  \bas
	\frac{1}{p} \io \big(\Teps(\cdot,T)+1\big)^p 
	\le \frac{1}{p} \cdot \bigg\{ \io \big(\Teps(\cdot,T)+1\big) \bigg\}^p \cdot |\Om|^{1-p}
	\qquad \mbox{for all $T>0$ and } \eps\in (0,1)
  \eas
  by the H\"older inequality, in view of (\ref{3.3}) we obtain (\ref{5.1}) from (\ref{5.6}).
\qed
Once more due to Lemma \ref{lem44}, this firstly implies estimates for the temperature variable
in some Lebesgue spaces involving conveniently large summability powers.
\begin{lem}\label{lem6}
  Let $q\in (0,3)$ and $T>0$. Then there exists $C(q,T)>0$ such that
  \be{6.1}
	\int_0^T \io (\Teps+1)^q
	\le C(q,T)
	\qquad \mbox{for all } \eps\in (0,1).
  \ee
\end{lem}
\proof
  Without loss of generality assuming that $q\in (2,3)$, we apply Lemma \ref{lem5} to $p:=q-2\in (0,1)$ to find
  $c_1=c_1(q,T)>0$ such that
  \be{6.2} 
	\int_0^T \io (\Teps+1)^{q-4} \Tepsx^2
	\le c_1
	\qquad \mbox{for all } \eps\in (0,1).
  \ee
  As Lemma \ref{lem44} yields $c_2=c_2(q)>0$ satisfying
  \bas
	\io (\Teps+1)^q
	\le c_2 \io (\Teps+1)^{q-4} \Tepsx^2
	+c_2
	\qquad \mbox{for all $t>0$ and } \eps\in (0,1),
  \eas
  we thus infer that
  \bas
	\int_0^T \io (\Teps+1)^q
	\le c_2 \int_0^T \io (\Teps+1)^{q-4} \Tepsx^2
	+ c_2 T
	\le c_1 c_2 + c_2 T
  \eas
  for all $T>0$ and $\eps\in (0,1)$.
\qed
When combined with Lemma \ref{lem6} in the course of another interpolation, the weighted gradient estimate in Lemma \ref{lem5} 
secondly implies bounds for the unweighted quantities $\Tepsx$.
\begin{lem}\label{lem7}
  If $r\in [1,\frac{3}{2})$ and $T>0$, then one can fix $C(r,T)>0$ in such a way that
  \be{7.1}
	\int_0^T \io |\Tepsx|^r
	\le C(r,T)
	\qquad \mbox{for all } \eps\in (0,1).
  \ee
\end{lem}
\proof
  Since the inequality $r<\frac{3}{2}$ ensures that $\frac{5r-6}{r}<1$, we can pick $p=p(r)\in (0,1)$ such that
  $p>\frac{5r-6}{r}$, and then obtain that $(5-p)r<\big(5-\frac{5r-6}{r}\big)\cdot r=6$,
  meaning that $q\equiv q(r):=\frac{(2-p)r}{2-r}$ satisfies
  \bas
	3-q
	= \frac{3(2-r)-(2-p)r}{2-r}
	= \frac{6-(5-p)r}{2-r}
	>0.
  \eas
  As thus $q<3$, besides employing Lemma \ref{lem5} we may therefore draw on Lemma \ref{lem6} to see that given $T>0$
  we can find $c_1=c_1(r,T)>0$ and $c_2=c_2(r,T)>0$ such that
  \be{7.2}
	\int_0^T \io (\Teps+1)^{p-2} \Tepsx^2
	\le c_1
	\quad \mbox{and} \quad
	\int_0^T \io (\Teps+1)^q \le c_2
	\qquad \mbox{for all } \eps\in (0,1).
  \ee
  Since the fact that $r<2$ enables us to use Young's inequality to see that
  \bas
	\int_0^T \io |\Tepsx|^r
	&=& \int_0^T \io \Big\{ (\Teps+1)^{p-2} \Tepsx^2 \Big\}^\frac{r}{2} \cdot (\Teps+1)^\frac{(2-p)r}{2} \\
	&\le& \int_0^T \io (\Teps+1)^{p-2} \Tepsx^2
	+ \int_0^T \io (\Teps+1)^\frac{(2-p)r}{2-r}
	\qquad \mbox{for all } \eps\in (0,1),
  \eas
  in line with our definition of $q$ we conclude (\ref{7.1}) from (\ref{7.2}).
\qed
Once more explicitly relying on the inequality $\al<\frac{3}{2}$ in (\ref{f}), we can now use Lemma \ref{lem6}
to control the integral on the right of (\ref{02.1}), and thereby obtain a space-time estimate for the viscosity-related
dissipation expressed therein.
\begin{lem}\label{lem8}
  For all $T>0$ there exists $C(T)>0$ such that
  \be{8.1}
	\int_0^T \io \vepsx^2
	\le C(T)
	\qquad \mbox{for all } \eps\in (0,1).
  \ee
\end{lem}
\proof
  We once more return to Lemma \ref{lem02}, but this time we use (\ref{gaeps}) here to see that
  \bea{8.2}
	\frac{1}{2} \frac{d}{dt} \io \veps^2
	+ \frac{a}{2} \frac{d}{dt} \io \uepsx^2
	&=& - \eps \io \vepsxx^2
	- \eps a \io \uepsxx^2
	- \io \gaeps(\Teps) \vepsx^2
	+ \io \feps(\Teps) \vepsx \nn\\
	&\le& - k_\gamma \io \vepsx^2
	+ \io \feps(\Teps) \vepsx 
	\qquad \mbox{for all $t>0$ and } \eps\in (0,1).
  \eea
  Since Young's inequality and (\ref{feps}) ensure that
  \bas
	\io \feps(\Teps) \vepsx
	&\le& \frac{k_\gamma}{2} \io \vepsx^2
	+ \frac{1}{2k_\gamma} \io f^2(\Teps) \\
	&\le&° \frac{k_\gamma}{2} \io \vepsx^2
	+ \frac{K_f^2}{2k_\gamma} \io (\Teps+1)^{2\al}
	\qquad \mbox{for all $t>0$ and } \eps\in (0,1),
  \eas
  we thereby obtain that
  \bas
	\frac{d}{dt} \bigg\{ \io \veps^2 + a \io \uepsx^2 \bigg\}
	+ k_\gamma \io \vepsx^2
	\le \frac{K_f^2}{k_\gamma} \io (\Teps+1)^{2\al}
	\qquad \mbox{for all $t>0$ and } \eps\in (0,1)
  \eas
  and hence
  \bas
	k_\gamma \int_0^T \io \vepsx^2
	\le \io v_{0\eps}^2 + a \io u_{0\eps x}^2
	+ \frac{K_f^2}{k_\gamma} \int_0^T \io (\Teps+1)^{2\al}
	\qquad \mbox{for all $T>0$ and } \eps\in (0,1).
  \eas
  Since $\sup_{\eps\in (0,1)} \big\{ \io v_{0\eps}^2 + a \io u_{0\eps x}^2 \big\} <\infty$ by (\ref{ie}),
  and since the inequality $2\al<3$ implies that also 
  $\sup_{\eps\in (0,1)} \int_0^T \io (\Teps+1)^{2\al}<\infty$ according to Lemma \ref{lem6},
  this entails (\ref{8.1}).
\qed
In preparation for an Aubin-Lions type argument, let us also record some regularity properties of the time derivatives
appearing in (\ref{0eps}). The first two of these are covered by the following.
\begin{lem}\label{lem9}
  Whenever $T>0$, one can find $C(T)>0$ such that
  \be{9.1}
	\int_0^T \big\| \vepst(\cdot,t)\big\|_{(W_0^{2,2}(\Om))^\star}^2 dt 
	\le C(T)
	\qquad \mbox{for all } \eps\in (0,1)
  \ee
  and
  \be{9.2}
	\int_0^T \io \uepst^2 
	\le C(T)
	\qquad \mbox{for all } \eps\in (0,1).
  \ee
\end{lem}
\proof
  For fixed $\psi\in C_0^\infty(\Om)$ fulfilling $\io \psi^2 + \io \psi_x^2 + \io \psi_{xx}^2 \le 1$,
  from (\ref{0eps}) and the Cauchy-Schwarz inequality as well as (\ref{gaeps}) we obtain that
  \bas
	\bigg| \io \vepst \psi \bigg|
	&=& \bigg| - \eps\io \vepsxx \psi_{xx} 
	- \io \gaeps(\Teps) \vepsx \psi_x
	- a \io \uepsx \psi_x
	+ \io \feps(\Teps) \psi_x \bigg| \\
	&\le& \eps \|\vepsxx\|_{L^2(\Om)}
	+ K_\gamma \|\vepsx\|_{L^2(\Om)}	
	+ a\|\uepsx\|_{L^2(\Om)}
	+ \|\feps(\Teps)\|_{L^2(\Om)}
  \eas
  for all $t>0$ and $\eps\in (0,1)$,
  so that with some $c_1>0$, because of (\ref{feps}) we have
  \bas
	\int_0^T \big\| \vepst(\cdot,t) \big\|_{(W_0^{2,2}(\Om))^\star}^2 dt
	&\le& c_1 \eps^2 \int_0^T \io \vepsxx^2
	+ c_1 \int_0^T \io \vepsx^2
	+ c_1 \int_0^T \io \uepsx^2
	+ c_1 \int_0^T \io \feps^2(\Teps) \\
	&\le& c_1 \eps \int_0^T \io \vepsxx^2
	+ c_1 \int_0^T \io \vepsx^2
	+ c_1 \int_0^T \io \uepsx^2
	+ c_1 K_f^2 \int_0^T \io (\Teps+1)^{2\al}
  \eas
  for all $T>0$ and $\eps\in (0,1)$.
  Again since $2\al<3$, we may thus combine (\ref{3.4}) and (\ref{3.2}) with Lemma \ref{lem8} and Lemma \ref{lem6} to arrive
  at (\ref{9.1}).\abs
  The estimate in (\ref{9.2}) directly results from (\ref{3.5}) and (\ref{3.1}), because by (\ref{0eps}) and Young's inequality,
  \bas
	\uepst^2 
	\le 2\eps^2 \uepsxx^2 + 2\veps^2
	\le 2\eps\uepsxx^2 + 2\veps^2
	\qquad \mbox{in } \Om\times (0,\infty)
  \eas
  for all $\eps\in (0,1)$.
\qed
Also the time derivatives of the temperture variable enjoy $\eps$-independent bounds in some suitable large spaces.
\begin{lem}\label{lem10}
  Let $\lam>3$. Then for all $T>0$ there exists $C(\lam,T)>0$ such that
  \be{10.1}
	\int_0^T \big\| \Tepst(\cdot,t)\big\|_{(W^{1,\lam}(\Om))^\star} dt 
	\le C(\lam,T)
	\qquad \mbox{for all } \eps\in (0,1).
  \ee
\end{lem}
\proof
  Using that $W^{1,\lam}(\Om)\hra L^\infty(\Om)$, we can find $c_1=c_1(\lam)>0$ such that
  $\|\psi_x\|_{L^\lam(\Om)} + \|\psi\|_{L^\infty(\Om)} \le c_1$ for all $\psi\in C^1(\bom)$ 
  fulfilling $\|\psi\|_{W^{1,\lam}(\Om)} \le 1$.
  For any such $\psi$, the third equation in (\ref{0eps}) along with the H\"older inequality, (\ref{gaeps}), (\ref{feps})
  and Young's inequality implies that
  \bas
	\bigg| \io \Tepst \psi \bigg|
	&=& \bigg| - \io \Tepsx \psi_x
	+ \io \gaeps(\Teps) \vepsx^2 \psi
	- \io \feps(\Teps) \vepsx \psi \bigg| \\
	&\le& \|\Tepsx\|_{L^\frac{\lam}{\lam-1}(\Om)} \|\psi_x\|_{L^\lam(\Om)}
	+ K_\gamma \cdot \bigg\{ \io \vepsx^2 \bigg\} \cdot \|\psi\|_{L^\infty(\Om)} \\
	& & + K_f \cdot \bigg\{ \io (\Teps+1)^\al |\vepsx| \bigg\} \cdot \|\psi\|_{L^\infty(\Om)} \\
	&\le& c_1 \|\Tepsx\|_{L^\frac{\lam}{\lam-1}(\Om)}
	+ c_1 K_\gamma \io \vepsx^2
	+ c_1 K_f \io (\Teps+1)^\al |\vepsx| \\
	&\le& c_1 \io |\Tepsx|^\frac{\lam}{\lam-1} + c_1 
	+ c_1 (K_\gamma+1) \io \vepsx^2
	+ \frac{c_1 K_f^2}{4} \io (\Teps+1)^{2\al}
  \eas
  for all $t>0$ and $\eps\in (0,1)$.
  Since the assumption $\lam>3$ warrants that $\frac{\lam}{\lam-1} \in (1,\frac{3}{2})$, and since $2\al<3$,
  an integration in time relying on Lemma \ref{lem7}, Lemma \ref{lem8} and Lemma \ref{lem6} yields the claim.
\qed
Now by means of an essentially straightforward subsequence extraction, we can indeed identify some limit functions
$v,u$ and $\Theta$ which satisfy the first sub-problem of (\ref{0}) in the sense specified in Definition \ref{dw}.
\begin{lem}\label{lem11}
  There exist $(\eps_j)_{j\in\N} \subset (0,1)$ and functions
  \be{11.1}
	\lbal
	v\in L^\infty((0,\infty);L^2(\Om)) \cap L^2_{loc}([0,\infty);W_0^{1,2}(\Om)), \\[1mm]
	u\in C^0(\bom\times [0,\infty)) \cap L^\infty((0,\infty);W_0^{1,2}(\Om))
	\qquad \mbox{and} \qquad \\[1mm]
	\Theta\in L^\infty((0,\infty);L^1(\Om)) \cap \bigcap_{q\in [1,3)} L^q_{loc}(\bom\times [0,\infty))
		\cap \bigcap_{r\in [1,\frac{3}{2})} L^r_{loc}([0,\infty);W^{1,r}(\Om))
	\ear
  \ee
  such that $u(\cdot,0)=u_0$ and $\Theta\ge 0$ a.e.~in $\Om\times (0,\infty)$, 
  that
  \be{11.01}
	f(\Theta) \in L^2_{loc}(\bom\times [0,\infty)),
  \ee
  that $\eps_j\searrow 0$ as $j\to\infty$, and that
  \begin{eqnarray}
	& & \veps\to v
	\qquad \mbox{a.e.~in $\Om\times (0,\infty)$ and in } L^2_{loc}(\bom\times [0,\infty)),
	\label{11.2} \\
	& & \vepsx \wto v_x
	\qquad \mbox{in } L^2_{loc}(\bom\times [0,\infty)),
	\label{11.3} \\
	& & \ueps \to u
	\qquad \mbox{in } L^\infty_{loc}(\bom\times [0,\infty)),
	\label{11.4} \\
	& & \uepsx \wto u_x
	\qquad \mbox{in } L^2_{loc}(\bom\times [0,\infty)),
	\label{11.5} \\
	& & \Teps \to \Theta
	\qquad \mbox{a.e.~in $\Om\times (0,\infty)$ and in } L^q_{loc}(\bom\times [0,\infty))
	\mbox{ for all $q\in [1,3)$,}
	\label{11.6} \\
	& & \feps(\Teps) \to f(\Theta)
	\qquad \mbox{in } L^2_{loc}(\bom\times [0,\infty)),
	\label{11.66} \\
	& & \Tepsx \to \Theta_x
	\qquad \mbox{in } L^r_{loc}(\bom\times [0,\infty))
	\mbox{ for all $r\in [1,\frac{3}{2})$,}
	\label{11.7} 
  \end{eqnarray}
  as $\eps=\eps_j\searrow 0$.
  Moreover,
  \be{wv}
	- \int_0^\infty \io v\vp_t - \io u_{0t} \vp(\cdot,0)
	= - \int_0^\infty \io \gamma(\Theta) v_x \vp_x
	- a \int_0^\infty \io u_x \vp_x
	+ \int_0^\infty \io f(\Theta) \vp_x
  \ee
  for all $\vp\in C_0^\infty(\Om\times [0,\infty))$, and
  \be{11.8}
	u_t=v
	\qquad \mbox{a.e.~in } \Om\times (0,\infty).
  \ee
\end{lem}
\proof
  From Lemma \ref{lem8} and Lemma \ref{lem9} we know that
  \bas
	(\veps)_{\eps\in (0,1)} 
	\ \mbox{ is bounded in } L^2((0,T);W_0^{1,2}(\Om))
	\quad \mbox{for all } T>0,
  \eas
  and that
  \bas
	(\vepst)_{\eps\in (0,1)} 
	\ \mbox{ is bounded in } L^2\big((0,T);(W_0^{2,2}(\Om))^\star\big)
	\quad \mbox{for all } T>0,
  \eas
  while Lemma \ref{lem3} together with Lemma \ref{lem9} asserts that
  \bas
	(\ueps)_{\eps\in (0,1)} 
	\ \mbox{ is bounded in } L^\infty((0,T);W_0^{1,2}(\Om))
	\quad \mbox{for all } T>0,
  \eas
  and that
  \bas
	(\uepst)_{\eps\in (0,1)} 
	\ \mbox{ is bounded in } L^2(\Om\times (0,T))
	\quad \mbox{for all } T>0.
  \eas
  Since from Lemma \ref{lem7}, Lemma \ref{lem3} and Lemma \ref{lem10} it follows that, moreover,
  \bas
	(\Teps)_{\eps\in (0,1)} 
	\ \mbox{ is bounded in } L^r((0,T);W^{1,r}(\Om))
	\quad \mbox{for all $T>0$ and any $r\in (1,\frac{3}{2})$,}
  \eas
  and that
  \bas
	(\Tepst)_{\eps\in (0,1)} 
	\ \mbox{ is bounded in } L^1\big((0,T);(W^{1,\lam}(\Om))^\star\big)
	\quad \mbox{for all $T>0$ and each } \lam>3,
  \eas
  three applications of an Aubin-Lions lemma (\cite{temam}), relying on the compactness of the embeddings 
  $W_0^{1,2}(\Om) \hra L^2(\Om)$, $W_0^{1,2}(\Om) \hra C^0(\bom)$ and $W^{1,r}(\Om) \hra L^1(\Om)$ for $r>1$,
  readily yield $(\eps_j)_{j\in\N}\subset (0,1)$ as well as $v\in L^2_{loc}([0,\infty);W_0^{1,2}(\Om))$,
  $u\in C^0(\bom\times [0,\infty)) \cap L^\infty((0,\infty);W_0^{1,2}(\Om))$ and 
  $\Theta \in \bigcap_{r\in (1,\frac{3}{2})} L^r_{loc}([0,\infty);W^{1,r}(\Om))$
  such that $\eps_j\searrow 0$ as $j\to\infty$, that (\ref{11.2})-(\ref{11.5}) and (\ref{11.7}) hold 
  and moreover $\Teps\to\Theta$ a.e.~in $\Om\times (0,\infty)$ as $\eps=\eps_j\searrow 0$, 
  whence also $u(\cdot,0)=u_0$ in $\Om$ and $\Theta\ge 0$ a.e.~in $\Om\times (0,\infty)$.
  As, apart from that, Lemma \ref{lem3} and Lemma \ref{lem6} imply that $(\veps)_{\eps\in (0,1)}$ is bounded in 
  $L^\infty((0,\infty);L^2(\Om))$, and that $(\Teps)_{\eps\in (0,1)}$ is bounded in $L^\infty((0,\infty);L^1(\Om))$
  and in $L^q(\Om\times (0,T))$ for all $T>0$ and each $q\in (1,3)$, it follows that the obtained limit functions 
  actually satisfy (\ref{11.1}), and that thanks to the Vitali convergence theorem also (\ref{11.6}) holds.
  Furthermore, again thanks to Vitali's theorem, the pointwise convergence features in (\ref{11.6}) ensure that, by (\ref{gflimit}),
  \be{11.12}
	\gaeps(\Teps) \to \gamma(\Theta)
	\quad \mbox{in } L^2_{loc}(\bom\times [0,\infty))
	\qquad \mbox{as } \eps=\eps_j \searrow 0,
  \ee
  and that also (\ref{11.66}) and (\ref{11.01}) are valid, for if in line with our hypothesis $\al<\frac{3}{2}$
  we fix any $\mu>0$ fulfilling $\mu\al<3$, then 
  \bas
	\sup_{\eps\in (0,1)} \int_0^T \io |\feps(\Teps)|^\mu \le K_f^\mu \sup_{\eps\in (0,1)} \int_0^T \io (\Teps+1)^{\mu\al}
	<\infty
	\qquad \mbox{for all } T>0
  \eas
  due to (\ref{feps}) and Lemma \ref{lem6}.\abs
  The verification of (\ref{wv}) thereupon becomes straightforward:
  For fixed $\vp\in C_0^\infty(\Om\times [0,\infty))$, the first equation in (\ref{0eps}) implies that
  \bea{11.9}
	\hs{-8mm}
	- \int_0^\infty \io \veps\vp_t - \io v_{0\eps} \vp(\cdot,0) 
	&=& - \eps \int_0^\infty \io \veps \vp_{xxxx}
	- \int_0^\infty \io \gaeps(\Teps) \vepsx \vp_x \nn\\
	& & - a \int_0^\infty \io \uepsx \vp_x
	+ \int_0^\infty \io \feps(\Teps) \vp_x
	\qquad \mbox{for all } \eps\in (0,1),
  \eea
  where by (\ref{11.2}),
  \be{11.10}
	- \int_0^\infty \io \veps \vp_t \to - \int_0^\infty \io v\vp_t
	\quad \mbox{and} \quad
	- \eps \int_0^\infty \io \veps \vp_{xxxx}
	\to 0
	\qquad \mbox{as } \eps=\eps_j \searrow 0,
  \ee
  while (\ref{ie}) and (\ref{11.5}) guarantee that
  \be{11.11}
	- \io v_{0\eps} \vp(\cdot,0)
	\to - \io u_{0t} \vp(\cdot,0)
	\quad \mbox{and} \quad
	- a \int_0^\infty \io \uepsx\vp_x 
	\to \int_0^\infty \io u_x \vp_x
	\qquad \mbox{as } \eps=\eps_j \searrow 0,
  \ee
  respectively. Moreover, when combined with (\ref{11.3}) the properties (\ref{11.12}) and (\ref{11.66}) imply that
  \bas
	- \int_0^\infty \io \gaeps(\Teps) \vepsx \vp_x
	\to - \int_0^\infty \io  \gamma(\Theta) v_x \vp_x
	\quad \mbox{and} \quad
	\int_0^\infty \io \feps(\Teps) \vepsx \vp_x
	\to \int_0^\infty \io f(\Theta) v_x \vp_x
  \eas
  as $\eps=\eps_j\searrow 0$,
  meaning that, indeed, (\ref{11.9}) implies (\ref{wv}).\abs
  Since, finally, (\ref{0eps}) also entails that for arbitrary $\vp\in C_0^\infty(\bom\times [0,\infty))$ we have
  \bas
	- \int_0^\infty \io \ueps\vp_t = \eps \int_0^\infty \io \ueps \vp_{xx} + \int_0^\infty \io \veps \vp
	\qquad \mbox{for all } \eps\in (0,1),
  \eas
  using (\ref{11.4}) and (\ref{11.2}) we infer on letting $\eps=\eps_j\searrow 0$ that also (\ref{11.8}) is valid.
\qed
\mysection{Solution properties of $\Theta$. Proof of Theorem \ref{theo17}}\label{sect5}
In line with Lemma \ref{lem11}, for the proof of Theorem \ref{theo17} it is now sufficient
to make sure that also the second identity (\ref{wt}) in Definition \ref{dw}
is satisfied by the limit functions found above.
At its core, this will amount to an appropriate limit passage in the respective contributions related to the viscoelastic
heat production term $\gaeps(\Teps) \vepsx^2$ in the third equation in (\ref{0eps}),
and a key step toward this will consist in confirming that with $(\eps_j)_{j\in\N}$ as provided by Lemma \ref{lem11},
the weak convergence property in (\ref{11.3}) can be turned into a strong approximation feature of the form in (\ref{07}).\abs
This will be achieved by suitably exploiting the weak identity (\ref{wv}) that is already known to be valid 
at this point, and in order to cleanly prepare our testing procedures in this regard, we separately 
record some technical prerequisites related to certain Steklov averages of relevance here.
\begin{lem}\label{lem12}
  Let $(v,u,\Theta)$ be as in Lemma \ref{lem11}, and let
  \be{12.1}
	\whv(x,t):=\lball
	v(x,t),
	\qquad & x\in\Om, \ t>0, \\[1mm]
	u_{0t}(x),
	\qquad & x\in\Om, \ t<0, 
	\ear
  \ee
  as well as
  \be{12.2}
	\whu(x,t):=\lball
	u(x,t),
	\qquad & x\in\Om, \ t>0, \\[1mm]
	u_0(x)+tu_{0t}(x),
	\qquad & x\in\Om, \ t<0. 
	\ear
  \ee
  Then writing
  \be{12.3}
	(S_h \vp)(x,t):= \frac{1}{h} \int_{t-h}^t \vp(x,s) ds,
	\qquad x\in\Om, \, t\in\R, \, h>0, \, \vp\in L^1_{loc}(\bom\times\R),
  \ee
  we have
  \be{12.4}
	(S_h \whv_x)(x,t) = \frac{\whu(x,t)-\whu(x,t-h)}{h}
	\qquad \mbox{for a.e.~$(x,t)\in\Om\times \R$ and each } h>0,
  \ee
  and, moreover,
  \be{12.5}
	S_h \whv \wto v
	\qquad \mbox{in } L^2_{loc}(\bom\times [0,\infty))
  \ee
  and
  \be{12.6}
	S_h \whv_x \wto v_x
	\qquad \mbox{in } L^2_{loc}(\bom\times [0,\infty))
  \ee
  as $h\searrow 0$.
\end{lem}
\proof
 Due to a well-known property of Steklov averages 
  (\cite[Lemma 3.2]{dibenedetto}),
  from the inclusions $\whv\in L^2_{loc}(\bom\times\R)$ and $\whv_x\in L^2_{loc}(\bom\times\R)$, as implied by (\ref{12.1})
  and (\ref{11.1}), it follows that $S_h\whv \to \whv$ and $S_h \whv_x \to \whv_x$ in $L^2_{loc}(\bom\times\R)$ as $h\searrow 0$,
  and that thus, again by (\ref{12.1}) both (\ref{12.5}) and (\ref{12.6}) and valid.\abs
  To verify (\ref{12.4}), we first note that (\ref{12.2}) and (\ref{12.1}) ensure that
  \be{12.7}
	\whu_{xt}=\whv_x
	\qquad \mbox{a.e.~in } \Om\times\R,
  \ee
  because for arbitrary $\vp\in C_0^\infty(\R\times\R)$, (\ref{11.8}) together with the continuity of $u$ and the identity
  $u(\cdot,0)=u_0$ can readily be seen to entail that
  \bas
	\int_0^\infty \io u \vp_{xt} = - \int_0^\infty \io v\vp_x - \io u_0 \vp_x(\cdot,0),
  \eas
  so that by (\ref{12.2}) and three integrations by parts,
  \bas
	- \int_{-\infty}^\infty \io \whu_x \vp_t
	&=& \int_{-\infty}^\infty \io \whu \vp_{xt} \\
	&=& \int_0^\infty \io u \vp_{xt} + \int_{-\infty}^0 \io (u_0+tu_{0t}) \vp_{xt} \\
	&=& - \int_0^\infty \io v\vp_x - \io u_0 \vp_x(\cdot,0) 
	+ \int_{-\infty}^0 \io (u_0+tu_{0t}) \vp_{xt} \\
	&=& - \int_0^\infty \io v\vp_x - \io u_0 \vp_x(\cdot,0) 
	- \int_{-\infty}^0 \io u_{0t} \vp_x
	+ \io u_0 \vp_x(\cdot,0) \\
	&=& - \int_{-\infty}^\infty \io \whv \vp_x
	= \int_{-\infty}^\infty \io \whv_x \vp
  \eas
  thanks to (\ref{12.1}).\abs
  Now having (\ref{12.7}) at hand, given $h>0$ and $\vp\in C_0^\infty(\Om\times\R)$ we abbreviate 
  $\vp^{(h)}(x,t):=\int_t^{t+h} \vp(x,s) ds$ for $x\in\Om$ and $t\in\R$, and noting that $\vp^{(h)} \in C_0^\infty(\Om\times\R)$
  we then obtain using the Fubini theorem and drawing on the definition of the distributional derivative $(\whu_x)_t$ that
  \bas
	\int_{-\infty}^\infty \io \vp(x,t) \cdot \bigg\{ \int_{t-h}^t \whv_x(x,s) ds\bigg\} dxdt
	&=& \int_{-\infty}^\infty \io \vp(x,t) \cdot \bigg\{ \int_{t-h}^t \whu_{xt}(x,s) ds\bigg\} dxdt \\
	&=& \int_{-\infty}^\infty \io \bigg\{ \int_s^{s+h} \vp(x,t) dt \bigg\} \cdot \whu_{xt}(x,s) dxds \\
	&=& \int_{-\infty}^\infty \io \vp^{(h)}(x,s) \whu_{xt}(x,s) dxds \\
	&=& - \int_{-\infty}^\infty \io \vp^{(h)}_t(x,s) \whu_x(x,s) dxds \\
	&=& - \int_{-\infty}^\infty \io \vp(x,s+h)\whu_x(x,s) dxds \\
	& & + \int_{-\infty}^\infty \io \vp(x,s) \whu_x(x,s) dxds.
  \eas
  As, by a linear substitution,
  \bas
	- \int_{-\infty}^\infty \io \vp(x,s+h) \whu_x(x,s) dxds
	= - \int_{-\infty}^\infty \io \vp(x,s) \whu_x(x,s-h) dxds,
  \eas
  after a division by $h$ this shows that
  \bas
	\int_{-\infty}^\infty \io \vp(x,t) \cdot (S_h \whv_x)(x,t) dxdt
	= \int_{-\infty}^\infty \io \vp(x,s) \cdot \frac{\whu_x(x,s)-\whu_x(x,s-h)}{h} dxds
  \eas
  and thereby establishes (\ref{12.4}), as $\vp\in C_0^\infty(\Om\times\R)$ was arbitrary.
\qed
On the basis of this, we can derive the following from (\ref{wu}) by simply using Young's inequality:
\begin{lem}\label{lem13}
  Let $(v,u,\Theta)$, $\whv$ and $S_h$ be as in Lemma \ref{lem11} and Lemma \ref{lem12}, and let 
  $\zeta\in C_0^\infty([0,\infty))$ be nonincreasing and such that $\zeta(0)=1$.
  Then
  \bea{13.1}
	\int_0^\infty \io \zeta(t) \gamma(\Theta(x,t)) v_x^2(x,t) dxdt
	&\ge& \frac{1}{2} \int_0^\infty \io \zeta'(t) v^2(x,t) dxdt
	+ \frac{1}{2} \io u_{0t}^2(x) dx \nn\\
	& & + \int_0^\infty \io \zeta(t) f(\Theta(x,t)) v_x(x,t) dxdt  \nn\\
	& & - a \cdot \limsup_{h\searrow 0} \int_0^\infty \io \zeta(t) u_x(x,t) \cdot (S_h\whv_x)(x,t) dxdt.
  \eea
\end{lem}
\proof
  We first note that on the right-hand side of (\ref{wv}) we have
  \be{13.02}
	\gamma(\Theta) v_x \in L^2_{loc}(\bom\times [0,\infty))
	\qquad \mbox{and} \qquad
	u_x \in L^2_{loc}(\bom\times [0,\infty))
  \ee
  due to (\ref{11.1}) and (\ref{g}).
  As furthermore also $f(\Theta)\in L^2_{loc}(\bom\times [0,\infty))$ by (\ref{11.01}),
  according to a standard approximation argument we thus infer that (\ref{wv}) continues to hold actually for each
  $\vp\in L^2((0,\infty);W_0^{1,2}(\Om))$ which is such that $\vp_t\in L^2(\Om\times (0,\infty))$
  and that $\vp=0$ a.e.~in $\Om\times (T,\infty)$ with some $T>0$, whence from (\ref{11.1}) it follows that for arbitrary $h>0$
  we may take 
  \bas
	\vp(x,t):=\zeta(t)\cdot (S_h \whv)(x,t),
	\qquad (x,t)\in\Om\times (0,\infty),
  \eas
  as a test function here.
  Since (\ref{12.3}) implies that
  \bas
	\vp_t(x,t)=\zeta'(t) \cdot (S_h v)(x,t)
	+ \zeta(t) \cdot\frac{\whv(x,t)-\whv(x,t-h)}{h}
  \eas
  and
  \bas
	\vp_x(x,t)=\zeta(t)\cdot (S_h \whv_x)(x,t)
  \eas
  for a.e.~$(x,t)\in \Om\times (0,\infty)$, we thereby achieve the identity
  \bea{13.3}
	& & \hs{-20mm}
	- \int_0^\infty \io \zeta'(t) v(x,t) \cdot (S_h \whv)(x,t) dxdt \nn\\
	& & \hs{-16mm}
	- \frac{1}{h} \int_0^\infty \io \zeta(t) v(x,t) \whv(x,t) dxdt
	+ \frac{1}{h} \int_0^\infty \io \zeta(t) v(x,t) \whv(x,t-h) dxdt
	- \io u_{0t}^2(x) dx \nn\\
	&=& - \int_0^\infty \io \zeta(t) \gamma(\Theta(x,t)) v_x(x,t) \cdot (S_h \whv_x)(x,t) dxdt \nn\\
	& & - a \int_0^\infty \io  \zeta(t) u_x(x,t) \cdot (S_h \whv_x)(x,t) dxdt \nn\\
	& & + \int_0^\infty \io \zeta(t) f(\Theta(x,t)) \cdot (S_h\whv_x)(x,t) dxdt,
  \eea
  because $\zeta(0)=1$ and $(S_h\whv)(x,0)=\frac{1}{h} \int_{-h}^0 u_{0t}(x) dt=u_{0t}(x)$ for a.e.~$x\in\Om$ by (\ref{12.1}).
  Here, (\ref{12.5}) together with the inclusion $v\in L^2_{loc}(\bom\times [0,\infty))$ and the fact that
  $\zeta'\in C_0^\infty([0,\infty))$ implies that
  \be{13.4}
	- \int_0^\infty \io \zeta'(t) v(x,t) \cdot (S_h \whv)(x,t) dxdt
	\to - \int_0^\infty \io \zeta'(t) v^2(x,t) dxdt
	\qquad \mbox{as } h\searrow 0,
  \ee
  while from (\ref{13.02}), (\ref{11.01}) and (\ref{12.6}) we similarly obtain that
  \be{13.5}
	- \int_0^\infty \io \zeta(t) \gamma(\Theta(x,t)) v_x(x,t) \cdot (S_h \whv_x)(x,t) dxdt
	\to - \int_0^\infty \io \zeta(t) \gamma(\Theta(x,t)) v_x^2(x,t) dxdt
  \ee
  and
  \be{13.6}
	\int_0^\infty \io \zeta(t) f(\Theta(x,t)) \cdot (S_h\whv_x)(x,t) dxdt
	\to \int_0^\infty \io \zeta(t) f(\Theta(x,t)) v_x(x,t) dxdt
  \ee
  as $h\searrow 0$.
  Apart from that, on the left-hand side of (\ref{13.3}) we may once more rely on (\ref{12.1}) to see by means of a 
  linear substitution and Young's inequality that since $\zeta\ge 0$,
  \bas
	& & \hs{-12mm}
	- \frac{1}{h} \int_0^\infty \io \zeta(t) v(x,t) \whv(x,t) dxdt
	+ \frac{1}{h} \int_0^\infty \io \zeta(t) v(x,t) \whv(x,t-h) dxdt \nn\\
	&=& - \frac{1}{h} \int_0^\infty \io \zeta(t) v^2(x,t) dxdt
	+ \frac{1}{h} \int_0^\infty \io \zeta(t) v(x,t) \whv(x,t-h) dxdt \nn\\
	&\le& -\frac{1}{2h} \int_0^\infty \io \zeta(t) v^2(x,t) dxdt
	+ \frac{1}{2h} \int_0^\infty \io \zeta(t) \whv^2(x,t-h) dxdt \nn\\
	&=& - \frac{1}{2h} \int_0^\infty \io \zeta(t) v^2(x,t) dxdt
	+ \frac{1}{2h} \int_0^\infty \io \zeta(t+h) v^2(x,t) dxdt
	+ \frac{1}{2h} \int_{-h}^0 \io \zeta(t+h) u_{0t}^2(x) dxdt \nn\\
	&=& \frac{1}{2} \int_0^\infty \io \frac{\zeta(t+h)-\zeta(t)}{h} v^2(x,t) dxdt
	+ \frac{1}{2} \cdot \bigg\{ \int_{-h}^0 \zeta(t+h) dt \bigg\} \cdot \io u_{0t}^2(x) dx,
  \eas
  so that since clearly $\frac{\zeta(\cdot+h)-\zeta(\cdot)}{h} \wsto \zeta'$ in $L^\infty((0,\infty))$ and 
  $\int_{-h}^0 \zeta(t+h) dt \to \zeta(0)=1$ as $h\searrow 0$, we see that
  \bas
	& & \hs{-30mm}
	\limsup \bigg\{ - \frac{1}{h} \int_0^\infty \io \zeta(t) v(x,t)\whv(x,t) dxdt
	+ \frac{1}{h} \int_0^\infty \io \zeta(t) v(x,t) \whv(x,t-h) dxdt \bigg\} \nn\\
	&\le& \frac{1}{2} \int_0^\infty \io \zeta'(t) v^2(x,t) dxdt
	+ \frac{1}{2} \io u_{0t}^2(x)dx.
  \eas  
  In conjunction with (\ref{13.4})-(\ref{13.6}), this shows that (\ref{13.3}) entails (\ref{13.1}).
\qed
Now an estimate for the last summand in (\ref{13.1}) can be obtained by now using the properties in (\ref{12.4}) and (\ref{12.2}),
and again utilizing Young's inequality an an adequate manner.
\begin{lem}\label{lem14}
  Let $\zeta\in C_0^\infty([0,\infty))$ be such that $\zeta'\le 0$ and $\zeta(0)=1$.
  Then with $(v,u,\Theta)$, $\whv$ and $S_h$ taken from Lemma \ref{lem11} and Lemma \ref{lem12}, we have
  \be{14.1}
	\limsup_{h\searrow 0} \int_0^\infty \io \zeta(t) u_x(x,t) \cdot (S_h\whv_x)(x,t) dxdt
	\le - \frac{1}{2} \int_0^\infty \io \zeta'(t) u_x^2(x,t) dxdt
	- \frac{1}{2} \io u_{0x}^2(x) dx.
  \ee
\end{lem}
\proof
  According to (\ref{12.4}) and (\ref{12.2}), we can rewrite
  \bea{14.2}
	\int_0^\infty \io \zeta(t) u_x(x,t) \cdot (S_h \whv_x)(x,t) dxdt
	&=& \frac{1}{h} \int_0^\infty \io \zeta(t) \whu_x(x,t) \cdot \big\{ \whu_x(x,t)-\whu_x(x,t-h)\big\} dxdt \nn\\
	&=& \frac{1}{h} \int_0^\infty \io \zeta(t) \whu_x(x,t-h) \cdot \big\{ \whu_x(x,t)-\whu_x(x,t-h)\big\} dxdt \nn\\
	& & + \frac{1}{h} \int_0^\infty \io \zeta(t) \cdot \big\{ \whu_x(x,t)-\whu_x(x,t-h)\big\}^2 dxdt \nn\\
	&=& \frac{1}{h} \int_0^\infty \io \zeta(t) \whu_x(x,t-h) \cdot \big\{ \whu_x(x,t)-\whu_x(x,t-h)\big\} dxdt \nn\\
	& & + h \int_0^\infty \io \zeta(t)\cdot (S_h\whv_x)^2(x,t) dxdt
	\qquad \mbox{for all } h>0.
  \eea
  here, in the second to last integral we can proceed in a way similar to that in Lemma \ref{lem13} to see by Young's inequality
  and a substitution that
  \bea{14.3}
	& & \hs{-26mm}
	\frac{1}{h} \int_0^\infty \io \zeta(t) \whu_x(x,t-h) \cdot \big\{ \whu_x(x,t)-\whu_x(x,t-h)\big\} dxdt \nn\\
	&=& \frac{1}{h} \int_0^\infty \io \zeta(t) \whu_x(x,t-h)\whu_x(x,t) dxdt \nn\\
	& & - \frac{1}{h} \int_0^\infty \io \zeta(t) \whu_x^2(x,t-h) dxdt \nn\\
	&\le& \frac{1}{2h} \int_0^\infty \io \zeta(t) \whu_x^2(x,t) dxdt
	- \frac{1}{2h} \int_0^\infty \io \zeta(t) \whu_x^2(x,t-h) dxdt \nn\\
	&=& \frac{1}{2h} \int_0^\infty \io \frac{\zeta(t)-\zeta(t+h)}{h} \whu_x^2(x,t) dxdt
	- \frac{1}{2} \cdot \bigg\{ \int_0^h \zeta(t) dt\bigg\} \cdot \io u_{0x}^2(x) dx \nn\\
	&\to& - \frac{1}{2} \int_0^\infty \io \zeta'(t) \whu_x^2(x,t) dxdt
	-\frac{1}{2} \io u_{0x}^2(x) dx 
	\qquad \mbox{as } h\searrow 0.
  \eea
  Since (\ref{12.6}) particularly entails that
  \bas
	\limsup_{h\searrow 0} \int_0^T \io (S_h\whv_x)^2(x,t) dxdt <\infty
	\qquad \mbox{for all } T>0,
  \eas
  however, the rightmost summand in (\ref{14.2}) satisfies
  \bas
	h \int_0^\infty \io \zeta(t)\cdot (S_h\whv_x)^2(x,t) dxdt
	\to 0
	\qquad \mbox{as } h\searrow 0,
  \eas
  so that (\ref{14.1}) becomes a consequence of (\ref{14.2}) and (\ref{14.3}).
\qed
In fact, the inequality announced in (\ref{06}) is implied by the previous two lemmata.
\begin{lem}\label{lem15}
  If $(v,u,\Theta)$ is as in Lemma \ref{lem11}, and if
  $\zeta\in C_0^\infty([0,\infty))$ is nonincreasing and such that $\zeta(0)=1$, then
  \bea{15.1}
	& & \hs{-20mm}
	\int_0^\infty \io \zeta(t) \gamma(\Theta(x,t)) v_x^2(x,t) dxdt
	- \frac{a}{2} \int_0^\infty \io \zeta'(t) u_x^2(x,t) dxdt \nn\\
	&\ge& \frac{1}{2} \io u_{0t}^2(x) dx
	+ \frac{a}{2} \io u_{0x}^2(x) dx \nn\\
	& & + \frac{1}{2} \int_0^\infty \io \zeta'(t) v^2(x,t) dxdt
	+ \int_0^\infty \io \zeta(t) f(\Theta(x,t)) v_x(x,t) dxdt.
  \eea
\end{lem}
\proof
  We only need to combine Lemma \ref{lem13} with Lemma \ref{lem14}.
\qed
The main technical objective of this section can now be achieved by suitably using (\ref{15.1}) in conjunction with
basic lower semicontinuity properties of weak convergence:
\begin{lem}\label{lem16}
  Let $(v,u,\Theta)$ and $(\eps_j)_{j\in\N}$ be as in Lemma \ref{lem11}, and let $T>0$.
  Then
  \be{16.1}
	\sqrt{\gaeps(\Teps)} \vepsx \to \sqrt{\gamma(\Theta)} v_x
	\quad \mbox{in } L^2(\Om\times (0,T))
	\qquad \mbox{as } \eps=\eps_j\searrow 0.
  \ee
\end{lem}
\proof
  Given $T>0$, we fix any nonincreasing $\zeta\in C_0^\infty([0,\infty))$ such that $\zeta\equiv 1$ on $(0,T)$,
  and from (\ref{11.6}), (\ref{gaeps}), (\ref{gflimit}), (\ref{11.3}) and the boundedness of $\zeta$ and of $\supp \zeta$
  we then readily obtain that
  \be{16.002}
	\sqrt{\zeta \gaeps(\Teps)} \vepsx
	\wto \sqrt{\zeta \gamma(\Theta)} v_x
	\quad \mbox{in } L^2(\Om\times (0,T))
	\qquad \mbox{as } \eps=\eps_j\searrow 0,
  \ee
  and that thus, by lower semicontinuity of the norm in $L^2(\Om\times (0,\infty))$ with respect to weak convergence, for
  \be{16.02}
	I_{1,\eps}:=\int_0^\infty \io \zeta(t) \gaeps(\Teps(x,t)) \vepsx^2(x,t) dxdt,
	\qquad \eps\in (0,1),
  \ee
  we have
  \be{16.2}
	\liminf_{\eps=\eps_j\searrow 0} I_{1,\eps} 
	\ge I_1:= \int_0^\infty \io \zeta(t) \gamma(\Theta(x,t)) v_x^2(x,t) dxdt.
  \ee
  Similarly, using (\ref{11.5}) along with the nonpositivity of $\zeta'$ we see that
  \be{16.03}
	I_{2,\eps} := - \frac{a}{2} \int_0^\infty \io \zeta'(t) \uepsx^2(x,t) dxdt,
	\qquad \eps\in (0,1),
  \ee
  satisfies
  \be{16.3}
	\liminf_{\eps=\eps_j\searrow 0} I_{2,\eps}
	\ge I_2:= - \frac{a}{2} \int_0^\infty \io \zeta'(t) u_x^2(x,t) dxdt.
  \ee
  To make appropriate use of this, we go back to Lemma \ref{lem02} and see on multiplying (\ref{02.1}) by $\zeta$ that
  \bas
	& & \hs{-30mm}
	\frac{1}{2} \int_0^\infty \zeta(t) \cdot \bigg\{ \frac{d}{dt} \io \veps^2(x,t) dx \bigg\} dt
	+ \frac{a}{2} \int_0^\infty \zeta(t) \cdot \bigg\{ \frac{d}{dt} \io \uepsx^2(x,t)dx \bigg\} dt \nn\\
	& & \hs{-24mm}
	+ \int_0^\infty \io \zeta(t) \gaeps(\Teps(x,t)) \vepsx^2(x,t) dxdt \nn\\
	& & \hs{-24mm}
	+ \eps \int_0^\infty \io \zeta(t) \vepsxx^2(x,t) dxdt
	+ \eps a \int_0^\infty \io \zeta(t) \uepsxx^2(x,t) dxdt \nn\\
	&=& \int_0^\infty \io \zeta(t) \feps(\Teps(x,t)) \vepsx(x,t) dxdt
	\qquad \mbox{for all } \eps\in (0,1),
  \eas
  where two integrations by parts show that since $\zeta(0)=1$,
  \bas
	& & \hs{-30mm}
	\frac{1}{2} \int_0^\infty \zeta(t)\cdot\bigg\{ \frac{d}{dt} \io \veps^2(x,t) dx\bigg\} dt \\
	&=& - \frac{1}{2} \int_0^\infty \io \zeta'(t) \veps^2(x,t) dxdt
	- \frac{1}{2} \io v_{0\eps}^2(x) dx
	\qquad \mbox{for all } \eps\in (0,1),
  \eas
  and that
  \bas
	& & \hs{-30mm}
	\frac{a}{2} \int_0^\infty \zeta(t)\cdot\bigg\{ \frac{d}{dt} \io \uepsx^2(x,t) dx\bigg\} dt \\
	&=& - \frac{a}{2} \int_0^\infty \io \zeta'(t) \uepsx^2(x,t) dxdt
	- \frac{a}{2} \io u_{0\eps x}^2(x) dx
	\qquad \mbox{for all } \eps\in (0,1).
  \eas  
  In view of (\ref{16.02}) and (\ref{16.03}), we thus obtain that
  for all $\eps\in (0,1)$,
  \bea{16.4}
	I_{1,\eps} + I_{2,\eps}
	&=& \frac{1}{2} \io v_{0\eps}^2(x)dx 
	+ \frac{a}{2} \io u_{0\eps x}^2(x)dx \nn\\
	& & + \frac{1}{2} \int_0^\infty \io \zeta'(t) \veps^2(x,t) dxdt
	+ \int_0^\infty \io \zeta(t)\feps(\Teps(x,t)) \vepsx(x,t) dxdt \nn\\
	& & - \eps \int_0^\infty \io \zeta(t) \vepsxx^2(x,t) dxdt
	- \eps a \int_0^\infty \io \zeta(t) \uepsxx^2(x,t) dxdt,
  \eea
  where
  \bas
	\frac{1}{2} \io v_{0\eps}^2(x) dx
	\to \frac{1}{2} \io u_{0t}^2(x)dx
	\quad \mbox{and} \quad
	\frac{a}{2} \io u_{0\eps x}^2(x) dx
	\to \frac{a}{2} \io u_{0x}^2(x) dx
	\qquad \mbox{as } \eps=\eps_j\searrow 0
  \eas
  by (\ref{ie}), and where the boundedness of $\supp \zeta$ ensures that
  \bas
	\frac{1}{2} \int_0^\infty \io \zeta'(t) \veps^2(x,t) dxdt
	\to \frac{1}{2} \int_0^\infty \io \zeta'(t) v^2(x,t) dxdt
	\qquad \mbox{as } \eps=\eps_j\searrow 0
  \eas
  due to (\ref{11.2}), and that
  \bas
	\int_0^\infty \io \zeta(t) \feps(\Teps(x,t)) \vepsx(x,t) dxdt
	\to \int_0^\infty \io \zeta(t) f(\Theta(x,t)) v_x(x,t) dxdt
	\qquad \mbox{as } \eps=\eps_j\searrow 0
  \eas
  thanks to (\ref{11.3}) and (\ref{11.66}).
  From (\ref{16.4}) it therefore follows that
  \bas
	\limsup_{\eps=\eps_j\searrow 0} (I_{1,\eps} + I_{2,\eps})
	&\le& I_3 := \frac{1}{2} \io u_{0t}^2(x) dx + \frac{a}{2} \io u_{0x}^2(x)dx \\
	& & + \frac{1}{2} \int_0^\infty \io \zeta'(t) v^2(x,t) dxdt
	+ \int_0^\infty \io \zeta(t) f(\Theta(x,t)) v_x(x,t) dxdt,
  \eas
  whence an application of Lemma \ref{lem15} reveals that
  \be{16.5}
	\limsup_{\eps=\eps_j\searrow 0}
	( I_{1,\eps} + I_{2,\eps} )
	\le \int_0^\infty \io \zeta(t) \gamma(\Theta(x,t)) v_x^2(x,t) dxdt
	- \frac{a}{2} \int_0^\infty \io \zeta'(t) u_x^2(x,t) dxdt.
  \ee
  In light of (\ref{16.2}) and (\ref{16.3}), however, this is possible only if
  \be{16.6}
	I_{1,\eps} \to I_1
	\qquad \mbox{as } \eps=\eps_j\searrow 0,
  \ee
  for otherwise (\ref{16.2}) and (\ref{16.3}) would imply that we could find $c_1>0$ and a subsequence
  $(\eps_{j_i})_{i\in\N}$ of $(\eps_j)_{j\in\N}$ such that
  \bas
	I_{1,\eps} \ge I_1 + c_1
	\quad \mbox{and} \quad
	I_{2,\eps} \ge I_2 - \frac{c_1}{2}
	\qquad \mbox{for all } \eps\in (\eps_{j_i})_{i\in\N},
  \eas
  meaning that
  \bas
	I_{1,\eps} + I_{2,\eps} \ge I_1 + I_2 + \frac{c_1}{2}
	\qquad \mbox{for all } \eps\in (\eps_{j_i})_{i\in\N}
  \eas
  and thereby contradicting (\ref{16.5}).\abs
  But in line with the definitions of $(I_{1,\eps})_{\eps\in (0,1)}$ and $I_1$ in (\ref{16.02}) and (\ref{16.2}),
  the relation in (\ref{16.6}) can be combined with (\ref{16.002}) to find that, in fact,
  \bas
	\sqrt{\zeta \gaeps(\Teps)} \vepsx \to \sqrt{\zeta \gamma(\Theta)} v_x
	\quad \mbox{in } L^2(\Om\times (0,\infty))
	\qquad \mbox{as } \eps=\eps_j\searrow 0.
  \eas
  When restricted to $\Om\times (0,T)$, due to the fact that $\zeta\equiv 1$ in $(0,T)$ this particularly entails (\ref{16.1}).
\qed
It remains to document a rather straightforward verification of the fact that the latter implies that also (\ref{wt}) holds:\abs
\proofc of Theorem \ref{theo17}. \quad
  We take $(v,u,\Theta)$ as provided by Lemma \ref{lem11}, and then readily obtain from (\ref{11.1}), (\ref{11.01}), (\ref{11.8})
  and (\ref{g}) that not only (\ref{17.1}) and (\ref{17.2}), and hence especially (\ref{w1}) and (\ref{w2}),
  but also (\ref{w3}) holds.
  The property $u(\cdot,0)=u_0$ as well as the inequality $\Theta\ge 0$ have been asserted by Lemma \ref{lem11} already,
  and so has been the identity in (\ref{wu}) which, in fact, is equivalent to (\ref{wv}) due to (\ref{11.8}).\abs
  To complete the proof, it is thus sufficient to note that for fixed $\vp\in C_0^\infty(\bom\times [0,\infty))$
  an integration by parts in (\ref{0eps}) shows that
  \bea{17.99}
	- \int_0^\infty \io \Teps \vp_t
	- \io \Theta_{0\eps} \vp(\cdot,0)
	&=& - \int_0^\infty \io \Tepsx \vp_x 
	+ \int_0^\infty \io \gaeps(\Teps) \vepsx^2 \vp \nn\\
	& & - \int_0^\infty \io \feps(\Teps) \vepsx \vp
	\qquad \mbox{for all } \eps\in (0,1),
  \eea
  and that here, as $\eps=\eps_j\searrow 0$ we have
  \bas
	- \int_0^\infty \io \Teps \vp_t 
	\to - \int_0^\infty \io \Theta \vp_t
	\quad \mbox{and} \quad
	- \io \Theta_{0\eps} \vp(\cdot,0) \to - \io \Theta_0 \vp(\cdot,0)
  \eas
  by (\ref{11.6}) and (\ref{ie}), 
  \bas
	- \int_0^\infty \io \Tepsx \vp_x \to - \int_0^\infty \io \Theta_x \vp_x
	\quad \mbox{and} \quad
	- \int_0^\infty \io \feps(\Teps) \vepsx \vp \to - \int_0^\infty \io f(\Theta) u_{xt} \vp
  \eas
  by (\ref{11.7}), (\ref{11.66}), (\ref{11.3}) and (\ref{11.8}), and that, again by (\ref{11.8}),
  the outcome of Lemma \ref{lem16} asserts that furthermore
  \bas
	\int_0^\infty \io \gaeps(\Teps) \vepsx^2 \vp
	\to \int_0^\infty \io \gamma(\Theta) u_{xt}^2 \vp
  \eas
  as $\eps=\eps_j\searrow 0$.
  Therefore, namely, (\ref{wt}) results from (\ref{17.99}).
\qed

\bigskip

{\bf Acknowlegements.} \quad
The author acknowledges support of the Deutsche Forschungsgemeinschaft (Project No. 444955436).
He moreover declares that he has no relevant financial or non-financial interests to disclose.\abs
{\bf Data availability statement.} \quad
Data sharing is not applicable to this article as no datasets were
generated or analyzed during the current study.

\end{document}